\documentclass[3p]{elsarticle}

\usepackage[T1]{fontenc}
\usepackage{amsmath}
\usepackage{amssymb}
\usepackage{amsthm}
\usepackage{color}
\usepackage{graphicx}
\usepackage{multirow}
\usepackage{algpseudocode} 
\usepackage{algorithm} 
\usepackage[]{subfig}
\usepackage{accents}
\renewcommand{\vec}{\mathrm{vec}}
\newcommand{\z}{\phantom{0}}
\renewcommand{\d}{\mathrm{d}}
\newcommand{\vect}[1]{\boldsymbol{#1}}

\newcommand{\X}{\mathcal{X}}
\newcommand{\Y}{\mathcal{Y}}
\newcommand{\G}{\mathcal{G}}

\graphicspath{{./PLOT/}}

\newtheorem{prop}{Proposition}
\newtheorem{rmk}{Remark}

\numberwithin{equation}{section}

\renewcommand{\d}{\mathrm{d}}

\newcommand{\kron}{\otimes}

\newcommand{\ds}{\ d s}

\renewcommand{\S}{\mathcal{S}}

\newcommand{\matb}[1]{\mathbf{#1}}

\newcommand{\core}[1]{\mathfrak{#1}}
\newcommand{\factor}[1]{\mathbf{#1}}
\newcommand{\trunc}[1]{\mathsf{#1}}

\begin{document}
\begin{frontmatter}  

 \title {A low-rank isogeometric solver based on Tucker
   tensors}

\author[mat,imati]{M. Montardini} \ead{monica.montardini@unipv.it}
\author[mat,imati]{G. Sangalli\corref{cor1}} \ead{giancarlo.sangalli@unipv.it}
\author[mat,imati]{M. Tani} \ead{mattia.tani@unipv.it}

\address[mat]{Dipartimento di Matematica, Universit\`a degli Studi di
Pavia, via Ferrata 5, Pavia, Italy}

\address[imati]{Istituto di Matematica Applicata e Tecnologie
  Informatiche ``E. Magenes'' del CNR, via Ferrata 5/a, Pavia, Italy}

\cortext[cor1]{Corresponding author}
\
%

\begin{abstract}

We propose an isogeometric solver  for Poisson problems that combines i)
low-rank tensor techniques to approximate the unknown solution and the system matrix, as a sum
of a few terms having Kronecker product structure, ii)  a Truncated Preconditioned
Conjugate Gradient solver to keep the rank of the iterates low, and
iii) a novel  low-rank preconditioner, based on the Fast
Diagonalization method where the eigenvector multiplication is  approximated
 by the  Fast Fourier Transform. Although the proposed strategy is written in arbitrary
dimension, we focus on the three-dimensional case and adopt the
Tucker format for low-rank tensor representation, which is well
suited in low dimension. We show 	 by  numerical tests that this choice
guarantees significant memory saving  compared to the full
tensor representation. We also extend and test the proposed strategy to linear elasticity problems.

 \end{abstract}

\end{frontmatter} 
 
 \section{Introduction}

 About twenty years ago, Tom Hughes  envisioned a revolutionary
 program for redesigning numerical  simulation methods by leveraging
 the knowledge and tools developed in computer-aided design (CAD),
 with the ultimate goal of unifying these two fields.  Isogeometric
 Analysis (IgA) appeared  in the seminal  paper \cite{Hughes2005} and
 since then a  community emerged,  that
brought together complementary interests and expertises: 
specialists in Computational Solid and Fluid Mechanics have been exploiting the potential of dealing with complex
computational domains directly from their spline or NURBS
representation within CAD;  the
geometric design community has been supporting  this challenging goal by
developing new analysis-suitable spline-based geometry parametrizations;
numerical analysts have been regaining momentum in studying spline
properties and their use in numerical solvers.

Isogeometric approximation, in particular, benefits from the regularity of
splines (see \cite{Beirao2014,Cottrell2007,Sande2020}),  a fundamental feature in
geometric design. Furthermore, the
tensor-product construction, widely used in multivariate
parameterizations, allows for computationally efficient methods when
utilised in the formation and solution of isogeometric linear
systems (see, e.g., \cite{HST_at_kunoth2018isogeometric} and
references therein). This is the context for the present work, which addresses
the use of low-rank tensor techniques.

The first contribution in this direction has been  given in 
\cite{Mantzaflaris2017}: there the authors use a low-rank approximated
representation of  the
coefficients of the Galerkin matrix, that incorporates the effect of
the geometry parametrization. An alternative approach, purely
algebraic, can be found in \cite{Hofreither2018}, and it is based on a
low-rank approximation of a small and dense matrix containing all the
non-zero entries of the original Galerkin matrix.  In these works, low-rank approximations are exploited for the
formation of the linear system, in order to write the Galerkin matrix
as a sum of a few matrices in Kronecker form. In this regard, we also
mention \cite{Juttler2017,Pan2019}, where the authors propose a
strategy to construct the geometry mapping that aims at minimizing the
number of Kronecker terms in this sum. More recently, low-rank tensor methods have been proposed for the
solution of IgA linear system in \cite{Georgieva2019} and
\cite{Bunger2020}, based respectively on an alternating least square
solver and  tensors in  Tucker format, and  on an alternating minimal
energy  method with  tensor-trains approximation of the unknown.

In the present paper we are mainly interested the three-dimensional Poisson problem. However, in principle the proposed methods can be generalized to an arbitrary number of dimensions $d$, and therefore they are often presented in this general setting.
We design a  low-rank isogeometric solver  whose computational cost
depends on the rank of the approximation and grows (almost) linearly
with respect to $n$, the number of
degrees of freedom (dofs) in each  space direction. First,  we
approximate the Galerkin matrix and the right-hand side with the
techniques from \cite{Driscoll2014,Dolgov2021}, based on  low-rank
Chebychev polyniomials.  The linear system is then solved iteratively
by a  Truncated Preconditioned Conjugate Gradient (TPCG) method
\cite{Kressner2011,Simoncini2022}, an extension of the  standard Preconditioned Conjugate Gradient 
method (PCG) where the tensor rank of the iterates is truncated.
We also introduce a  preconditioner
inspired by the Fast Diagonalization method
\cite{Lynch1964,Sangalli2016} where, however,  the projection onto
the eigenvectors basis is computed by a   Fast Fourier Transform  and
the eigenvalues are replaced  by a suitable low-rank approximation
using the results of \cite{Braess2005,Hackbusch2019}. The proposed  preconditioner yields a setup/application cost which is (almost) linear with respect to $n$, and it is robust with respect to the mesh size and the spline degree.   We also extend the proposed low-rank solving strategy to linear elasticity problems. 

We choose the Tucker format to represent tensors, which is considered
well suited in low dimension \cite{Kressner2017,Dolgov2021}.  
We remark that for higher dimension other low-rank tensor formats, such as hierarchical Tucker or tensor trains,  are considered  more efficient.
In any case, our numerical experience confirms that the use of a low-rank Tucker approximation
offers advantages with respect to the full tensor approach. These
advantages are problem dependent but always significant, with order of magnitude
speedup and memory saving.

The preconditioner we propose is innovative as it leverages the
low-rank structure of the solver for optimal efficiency. However,
there are many approaches to preconditioning the problem in its
full-rank formulation, particularly methods that exploit tensorial
properties and are robust with respect to polynomial degree: we
refer to the works  \cite{Beirao2017,Donatelli2015,Garcia2018,Hofreither2017,Sangalli2016,Tielen2020}  and the references contained within them.

The paper is organized as follows. In Section \ref{sec:preliminaries} we present the basics of IgA and of tensor calculus. The model problem is introduced in Section \ref{sec:model_problem}. The core of the paper is Section \ref{sec:low_rank_solver}, where we recall the Truncated Preconditioned Conjugate Method, the truncation operators and where we introduce the novel preconditioning strategy. We present some numerical experiments in Section \ref{sec:numerics} while in the last section we draw some conclusions. We rely on the two appendices the  technical parts: the approximation of the linear system in Tucker format in \ref{sec:appendix2} and the approximation of the eigenvalues and eigenvectors needed for building the preconditioner in \ref{sec:appendix}.

\medskip

\emph{This paper is submitted to the CMAME Special Issue in Honor of the
Lifetime Achievements of  Dr. Thomas J.R. Hughes. The authors would like  to congratulate
Tom  for his extraordinary scientific
contributions and express their gratitude for his support in
shaping their work with his visionary ideas.
}

\section{Preliminaries}
\label{sec:preliminaries}
\subsection{B-Splines} 
\label{sec:bsplines}

A knot vector in $[0,1]$ is   a sequence of non-decreasing points
 $\Xi:=\left\{ 0=\xi_1 \leq \dots \leq \xi_{m+p+1}=1\right\}$, where   $m$ and $p$ are two positive integers that represent the number of basis functions associated to the knot vector and their polynomial degree, respectively.
 We focus on open knot vectors, i.e.  we set $\xi_1=\dots=\xi_{p+1}=0$
 and $\xi_{m}=\dots=\xi_{m+p+1}=1$. 
According to   Cox-De Boor recursion formulas (see  \cite{DeBoor2001}), univariate B-splines $\widehat{b}_{i}^{(p)}: [0,1]  \rightarrow \mathbb{R}$ are  piecewise polynomials defined  for $i=1,\dots,m$ as \\
\indent for $p=0$:
\begin{align*} 
\widehat{b}_{i}^{(0)}(\eta) = \begin{cases}1 &  { \textrm{if }} \xi_{i}\leq \eta<\xi_{i+1},\\
0 & \textrm{otherwise,}
\end{cases}
\end{align*}
\indent for $p\geq 1$: 

$$
 \widehat{b}_{i}^{(p)}(\eta)= \dfrac{\eta-\xi_{i}}{\xi_{i+p}-\xi_{i}}\widehat{b}_{  i}^{(p-1)}(\eta)   +\dfrac{\xi_{i+p+1}-\eta}{\xi_{i+p+1}-\xi_{i+1}}\widehat{b}^{(p-1)}_{  i+1}(\eta)   \quad { \textrm{if }} \xi_{i}\leq  \eta<\xi_{i+p+1},  
$$

 where we assume $0/0=0$. 
The univariate spline space is defined as
\begin{equation*}
\widehat{\S}_{h}^p(\Xi) : = \mathrm{span}\{\widehat{b}_{i}^{(p)}\}_{i = 1}^m,
\end{equation*}
where $ h:=\max\{  \xi_{i+1}-\xi_i  \ | \ i=1,\dots,m+p \}$  denotes the mesh-size.
The smoothness of the B-splines at the interior knots is determined by their multiplicity (see \cite{DeBoor2001}). We refer to \cite{Cottrell2009} for more properties on B-splines.
We write $\widehat{b}_i$ instead of $	\widehat{b}^{(p)}_{i}$, when the degree is clear from the context.
 Multivariate B-splines are defined as tensor product of univariate
B-splines. 
In particular, for $d$-dimensional problems, we introduce $d$ univariate knot vectors $\Xi_l:=\left\{ \xi_{l,1} \leq \dots \leq
  \xi_{l,m_l+p_l+1}\right\}$  for $l=1,\ldots, d$, where $m_l$ and $p_l$ are positive integers for $l=1,\dots,d$.
We denote with $h_{l}$ the mesh-size associated to the knot vector $\Xi_l$ for $l=1,\dots,d$, and we define by $h:=\max\{h_{l}\ | \ l=1,\dots, d\}$  the maximal mesh-size.
For simplicity, we assume $p_1=\dots=p_d=:p$, but the general case is similar.
We denote the $i$-th  univariate function in the $k$-th direction as $\widehat{b}_{k,i}$ for  $i=1,\dots,m_k$ and $k=1,\dots,d$.
We assume that all the knot vectors
are  uniform, i.e. the internal knots are equally spaced. 

The multivariate B-splines are defined as
\begin{equation*} 
\widehat{B}_{ \vect{i}}(\vect{\eta}) : =
\widehat{b}_{1,i_1}(\eta_1) \ldots \widehat{b}_{d,i_d}(\eta_d),
\end{equation*}
 where  $\vect{i}:=(i_1,\dots,i_d)$ and  $\vect{\eta} = (\eta_1, \ldots, \eta_d) \in \widehat{\Omega} := [0,1]^d$.   The  corresponding spline space  is defined as
\begin{equation*}
\widehat{\boldsymbol{\S}}^{p}_{h}  := \mathrm{span}\left\{\widehat{B}_{\vect{i}} \ \middle| \text{where } \vect{i}:=(i_1,\dots,i_d)\text{ and }\ i_l = 1,\dots, m_l \text{ for } l=1,\dots,d \right\}= {\widehat{\S}^{p}_{h_{1}}\otimes\dots\otimes\widehat{\S}^{p}_{h_{d}}},
\end{equation*}  
where $\widehat{\S}^{p}_{h_{k}} :=\text{span}\{\widehat{b}_{k,i}\ | \ i=1,\dots,m_k\}$ for $k=1,\dots,d$.

\subsection{Isogeometric spaces}
\label{sec:iso_space}
We assume that our computational domain $\Omega\subset\mathbb{R}^d$ is given by a spline parametrization  $\vect{F}\in \widehat{\boldsymbol{\S}}^{p}_{h} $, i.e. $\Omega = \vect{F}(\widehat{\Omega})$. We also assume that $\vect{F}$ has a non-singular Jacobian everywhere.
%
For our purpose we also need  the spline space with homogeneous Dirichlet boundary
 condition:
\begin{equation*}
\widehat{\boldsymbol{\S}}^p_{h,0} :=\left\{ \widehat{v}_h\in \widehat{\boldsymbol{\mathcal{S}}}^{{p}}_h  \ \middle| \ \widehat{v}_h = 0 \text{ on } \partial\widehat{\Omega}  \right\}.
\end{equation*}
By introducing a colexicographical reordering of  the basis functions, we can  write  

\begin{align}
\label{eq:spline_space}
  \widehat{\boldsymbol{\S}}^p_{h,0}  & = \ \text{span}\left\{ \widehat{b}_{1,i_1}\dots\widehat{b}_{d,i_d} \ \middle| \ i_l = 2,\dots , m_l-1; \ l=1,\dots,d\ \right\}  =\text{span}\left\{ \widehat{B}_{i} \ \middle|\ i =1,\dots , N_{dof}   \ \right\}\\
  & = {\widehat{\S}^p_{h_1,0} \otimes\dots\otimes\widehat{\S}^p_{h_d,0}},
\end{align}
where  $\widehat{\S}^p_{h_l,0} :=\text{span}\left\{\widehat{b}_{l,i} \ | \ i=2,\dots,m_l-1 \right\}$ for $l=1,\dots,d$,  $N_{dof}:=n_1\dots n_d$, $n_l:=m_l-2$, and, where  with a little abuse of notation, we identify the multi-index $\vect{i} = (i_1,\dots, i_d)$ with the scalar index $i = i_1-1+ \sum_{l=2}^d (i_l-2)\prod_{k=1}^{l-1} n_k$.


Finally, the isogeometric space we consider is the isoparametric push-forward of $ \widehat{\boldsymbol{\S}}^p_{h,0} $ through the geometric map $\vect{F}$, i.e. we take
\begin{equation}
\boldsymbol{\mathcal{V}}_{h} := \text{span}\left\{  B_{i}:=\widehat{B}_{i}\circ \vect{F}^{-1} \ \middle| \ i=1,\dots , N_{dof}   \right\}.
\label{eq:disc_space}
\end{equation}

  \subsection{Tensors calculus}
  \label{sec:tensor_calculus}
  We recall the essential properties of tensor calculus that are useful in our context. Further details can be found e.g. in the survey \cite{Kolda2009}.

  \paragraph{Tensors}
Given $d\in\mathbb{N}$, a \textit{tensor}  $\mathcal{X}\in \mathbb{R}^{n_1\times\dots\times n_d}$ with $n_1,\dots,n_d\in\mathbb{N}$ is a $d$-dimensional array. We denote its entries as   $[\mathcal{X}]_{i_1,\dots,i_d}\in\mathbb{R}$ for $i_l=1,\dots,n_l$ and $l=1,\dots,d$. Usually $d$ is called the \textit{order} of $\mathcal{X}$. Tensors of order 1 are vectors, while tensors of order 2 are matrices.

The scalar product of two tensors $\mathcal{X}, \mathcal{Y}\in \mathbb{R}^{n_1\times\dots\times n_d}$  is defined as
 $$<\mathcal{X}, \mathcal{Y}>:= \sum_{i_1=1}^{n_1}\dots\sum_{i_d=1}^{n_d}[\mathcal{X}]_{i_1,\dots,i_d}[\mathcal{Y}]_{i_1,\dots,i_d}.$$
 The associated norm is  called Frobenius norm and is denoted as
 $\|\cdot\|_F$.

 The {vectorization operator} ``vec'' applied to a tensor stacks its entries  into a column vector as
\[  [\text{vec}(\mathcal{X})]_{j}=[\mathcal{X}]_{i_1,\dots,i_{d}} \text{ for }   i_l=1,\dots,n_{l} \text{ and for } l=1,\dots,d,   \] 
where  $j:=i_1+\sum_{k=2}^{d}\left[(i_k-1)\Pi_{l=1}^{k-1}n_l\right]$.  Thanks to the ``vec'' operator, we can always represent a vector $\vect{x}\in\mathbb{R}^{n_1 \dots n_d}$  as a tensor $\mathcal{X}\in \mathbb{R}^{n_1\times\dots\times n_d}$ and viceversa:   \begin{equation}
\label{eq:vec_tens}
\vect{x}=\text{vec}(\mathcal{X}).
\end{equation}
Note that $\|\vect{x}\|_2=\|\mathcal{X}\|_F$, where $\|\cdot\|_2$ is the euclidean norm of vectors.
In the following, we denote the tensor associated to a vector with the calligraphic upper case version of the letter used to indicate the vector, e.g. the tensors $\mathcal{Y, X \text{ and } \G}$ correspond, respectively, to the vectors $\vect{y\text{, }x \text{ and }g}$. 

Assume for simplicity that $n_i = n$, for $i = 1,\ldots,d$. Then the number of entries in the tensor is $n^d$, which   grows exponentially with $d$.
This is the so-called "curse of dimensionality", which makes practically impossible to store explicitly a high order tensor. For this reason memory-efficient representations based on tensor products are preferable.

\paragraph{Kronecker product and $m$-mode product}
Given two matrices $\mathbf{C}\in\mathbb{R}^{n_1\times n_2}$ and $\mathbf{D}\in\mathbb{R}^{n_3\times n_4}$, their Kronecker product  is defined as
\begin{equation}
\label{eq:kron_matrices}\mathbf{C}\otimes \mathbf{D}:=\begin{bmatrix}
[\mathbf{C}]_{1,1}\mathbf{D}  & \dots& [\mathbf{C}]_{1,n_2}\mathbf{D}\\
\vdots& \ddots &\vdots\\
[\mathbf{C}]_{n_1, 1}\mathbf{D}& \dots & [\mathbf{C}]_{n_1, n_2}\mathbf{D}
\end{bmatrix}\in \mathbb{R}^{n_1n_3\times n_2 n_4},
\end{equation}
where   $[\mathbf{C}]_{i,j}$ denotes the $ij$-th entry of the matrix $\mathbf{C}$.  
  More generally, the Kronecker product of two tensors $\mathcal{C} \in \mathbb{R}^{m_1 \times \ldots \times m_d}$ and $\mathcal{D} \in \mathbb{R}^{n_1 \times \ldots \times n_d}$ is the tensor

\begin{equation}
\label{eq:kron_tensors}
\left[\mathcal{C}\otimes\mathcal{D}\right]_{k_1,\dots,k_d} = [\mathcal{C}]_{i_1,\dots,i_d} [\mathcal{D}]_{j_1,\dots,j_d}\in\mathbb{R}^{n_1m_1\times\dots\times n_dm_d}\quad \text{with}\quad k_l=j_l+(i_l-1) n_l \quad \text{for}\quad l=1,\dots,d.
\end{equation} 

 For $m=1,\dots,d$  we  introduce the {$m$-mode product}   of a tensor $\mathcal{X}\in\mathbb{R}^{n_1\times\dots\times n_{d}}$ with a matrix $\mathbf{J}\in\mathbb{R}^{\ell\times n_m}$,   that we denote by  $\mathcal{X}\times_m \mathbf{J}$.  This is 
  a tensor of size $n_1\times\dots\times n_{m-1}\times \ell \times n_{m+1}\times \dots n_{d}$,   whose entries are
\[\left[ \mathcal{X}\times_m \mathbf{J} \right]_{i_1, \dots, i_{d}} := \sum_{j=1}^{n_m} [\mathcal{X}]_{i_1,,\dots, i_{m-1},j,i_{m+1},\dots,i_{d}} [\mathbf{J}]_{j,i_m }.\]
The product between a vector and a Kronecker product matrix can be expressed in terms of a matrix-tensor product.
Precisely, given $\mathbf{J}_i\in\mathbb{R}^{\ell_i\times n_i}$ for $i=1,\dots, d$, it holds
\begin{equation}
\label{eq:kron_vec_multi}
\left(\mathbf{J}_{d}\otimes\dots\otimes \mathbf{J}_1\right)\text{vec}\left(\mathcal{X}\right)=\text{vec}\left(\mathcal{X} \times_d  \mathbf{J}_{d} \dots \times_1 \mathbf{J}_1\ \right).
\end{equation} 

\paragraph{Tucker format}
 
 In this paper, we use the Tucker format to represent tensors. We adopt the notations of \cite{Oseledets2009}.
 A tensor $\mathcal{X}\in \mathbb{R}^{n_1\times\dots\times n_d}$ is in Tucker format if it is expressed as
\begin{equation}
\label{eq:Tucker}
\mathcal{X}=\core{X}\times_d\factor{X}_d\dots\times_1\factor{X}_1.
\end{equation}
where $\core{X} \in \mathbb{R}^{r_1 \times \ldots \times r_d}$ is the core tensor and  $\factor{X}_k \in\mathbb{R}^{n_k\times r_k}$ are the factor matrices.
The $d-$uple $(r_1,\dots,r_d)$ is called the \textit{multilinear rank} of $\mathcal{X}$. The storage for a Tucker tensor is bounded by
$drn+r^d$ where, here and throughout the paper we denote
$$r:=\max_{i=1,\dots,d}{r_i}\qquad \text{and} \qquad  n:=\max_{i=1,\dots,d}{n_i}.$$ 
 
Thus, the Tucker format still suffers from an exponential increase in memory storage with respect to the dimension $d$. 
Nevertheless, if $r \ll n$, the memory storage of a Tucker tensor is much smaller than the one of a full tensor. 

  
One possible way to compute the Tucker format of  a tensor $\X$ is  to use the High Order Singular Value Decomposition (HOSVD) \cite{Tucker1966,De2000}, that is an extension to tensors of the  matrix  Singular Value Decomposition (SVD).   
An approximation of the  Tucker format of $\X$ can be found for
example by  the Sequentially Truncated HOSVD (ST-HOSVD) \cite{Vannieuwenhoven2012}.
The ST-HOSVD algorithm  computes a Tucker tensor  $\widetilde{\X}=\textsc{sthosvd}(\X,\epsilon)$ that approximates $\X$ and that  satisfies 
$$
\|\X-\widetilde{\X}\|_F\leq\epsilon\|\X\|_F
$$
for a given tolerance $\epsilon>0$.
We remark that the factor matrices resulting from the application of the ST-HOSVD are orthogonal.
 For other truncation algorithms, we refer to the surveys \cite{Kolda2009,Grasedyck2013}.

Due to  the relation between tensors and vectors provided by \eqref{eq:vec_tens}, we also say that a vector $\vect{x}\in \mathbb{R}^{n_1\dots n_d}$ is in Tucker format if it is expressed as 
\begin{align*}\vect{x}&= (\factor{X}_d \otimes\dots\otimes\factor{X}_1)\text{vec}(\core{X})\\
&=\sum_{i_d=1}^{r_d} \dots
                                                                                                \sum_{i_1=1}^{r_1}[\core{X}]_{i_1,\dots,i_d}
                                                                                                \vect{x}_d^{(i_d)}\otimes\dots\otimes
                                                                                                \vect 
                                                                                                x_1^{(i_1)},
\end{align*} 
where   $\vect{x}_k^{(i_k)}\in\mathbb{R}^{n_k}$   is the $i_k$-th column of $\factor{X}_k\in\mathbb{R}^{n_k\times r_k}$ and $\core{X}\in\mathbb{R}^{r_1\times\dots \times r_d}$. Note that the tensor $\mathcal{X}\in\mathbb{R}^{n_1\times\dots\times n_d}$ associated to $\vect{x}$  can thus be written as \eqref{eq:Tucker}. 
To simplify the exposition, we will sometimes refer to the multilinear rank of a vector in Tucker format, meaning with this the multilinear rank of the associated tensor.

Furthermore, we say that a matrix $\matb{C}\in\mathbb{R}^{(n_1\dots n_d)\times (n_1\dots n_d)}$ is in Tucker format if it is expressed as
 \begin{equation}
 \label{eq:tucker_matrix}
 \matb{{{C}}}=\sum_{i_d=1}^{R^C_{d}} \dots \sum_{i_1=1}^{R^C_{1}}[\core{C}]_{i_1,\dots,i_d} \matb{C}_{(d,i_d)}\otimes\dots\otimes \matb{C}_{(1,i_1)},
 \end{equation}
 where $\core{C}\in\mathbb{R}^{R^C_{1}\times\dots\times R^C_{d}}$ and $\matb{C}_{(k,i_k)}\in\mathbb{R}^{n_k\times n_k}$ for $i_k=1,\dots R^C_{k}$ and $k=1,\dots,d$.

 \paragraph{Binary operations for Tucker format}
 Let   $\matb{C}\in\mathbb{R}^{(n_1\dots n_d)\times (n_1\dots n_d)}$ be a matrix in Tucker format as in
 \eqref{eq:tucker_matrix} and $\vect{x}\in\mathbb{R}^{n_1\dots n_d}$ a Tucker vector represented by $\mathcal{X}$  as in \eqref{eq:Tucker}.  
The  multiplication between the matrix $\matb{C}$ and the vector $\vect{x}$ can be efficiently computed as
\begin{equation}
	\label{eq:tt_product}
\matb{C}\vect{x}=\text{vec}\left(\left(\core{C}\otimes\core{X}\right)\times_d  \left[\factor{C}_{(d,1)}\factor{X}_d, \dots , \factor{C}_{(d,R^C_d)}\factor{X}_d\right]\times_{d-1} \dots\times_1 \left[\factor{C}_{(1,1)}\factor{X}_1, \dots ,\factor{C}_{(1,R^C_1)}\factor{X}_1\right]  \right).
\end{equation}  
Note that the multilinear rank of $\matb{C}\vect{x}$  is $[R^C_{1}r_1,\dots,R^C_{d}r_d]$.

Let $\vect{y} \in\mathbb{R}^{n_1\dots n_d}$  be  a vector in Tucker format represented by $\mathcal{Y}:=\core{Y}\times_d \factor{Y}_d \dots \times_{1}\factor{Y}_1$ where $\core{Y}\in\mathbb{R}^{s_1\times\dots\times s_d}$ is the core tensor and $\factor{Y}_i\in\mathbb{R}^{n_i\times s_i}$ are the factor matrices. The scalar product between $\vect{x}$ and $\vect{y}$ can be computed exploiting their Tucker format as
$$
\vect{x}\cdot\vect{y}=<\mathcal{X},\mathcal{Y}> = \text{vec}(\core{X})^T(\factor{X}_d^T\factor{Y}_d\otimes \dots\otimes \factor{X}_1^T\factor{Y}_1)\text{vec}(\core{Y})=\text{vec}(\core{X})^T\text{vec}(\core{Y}\times_d\factor{X}_d^T\factor{Y}_d  \dots\times_1 \factor{X}_1^T\factor{Y}_1).
$$

The sum of  two Tucker vectors $\vect{x}$ and $\vect{y}$ is a Tucker vector $\vect{z}:=\vect{x}+\vect{y}$ represented by the Tucker tensor $\mathcal{Z}:=\core{Z}\times_d \left[\factor{X}_d, \factor{Y}_d\right]\times_{d-1}\dots \times_1\left[\factor{X}_1,  \factor{Y}_1\right]$, where the core tensor $\core{Z}\in\mathbb{R}^{(r_1+s_1)\times\dots\times (r_d+s_d)}$ is a block-diagonal tensor defined by concatenating on the diagonal the core tensors of   $\mathcal{X}$ and $\mathcal{Y}$. Note that the multilinear rank of $\mathcal{Z}$ is $(r_1+s_1, \dots, r_d+s_d)$, i.e. the sum of the multilinear ranks of the addends.

\section{Model problem}
 \label{sec:model_problem}

Our model problem is  the Poisson problem
\begin{equation*}
\left\{
	\begin{array}{rcll}
-\Delta u & = & f  & \text{ in } \Omega\\
u & = & 0 & \text{ on }\partial \Omega
	\end{array}
\right.
\end{equation*}
where $\Omega\subset\mathbb{R}^d$.   For the sake of simplicity, we
consider homogeneous Dirichlet boundary condition. The weak formulation reads: find $u\in H^1_0(\Omega)$ such that for all $v\in H^1_0(\Omega)$ it holds
$$
a(u,v)=F(v)
$$
where 
\begin{equation*}
a(u,v):=\int_{\Omega}\nabla u \cdot \nabla v\ \d\Omega \quad \text{and}\quad F(v):=\int_{\Omega}fv\ \d\Omega.
\end{equation*}
The isogeometric  discretization with the space \eqref{eq:disc_space} yields to the following discrete problem: find $u_h\in\boldsymbol{\mathcal{V}}_h$ such that for all $v_h\in\boldsymbol{\mathcal{V}}_h$ it holds
\begin{equation}
\label{eq:dicrete_bilinear}
a(u_h,v_h)= F(v_h).
\end{equation}
The linear system associated to \eqref{eq:dicrete_bilinear} is
\begin{equation}
\label{eq:linear_sys}
{\matb{A}}\vect{x}= {\vect{f}}
\end{equation}
where $[{\matb{A}}]_{i,j}:=a(B_i,B_j) $ and $[\vect{f}]_i:=F(B_i)$ for $i, j=1,\dots,N_{dof}$. 
 We approximate the system matrix with a matrix in Tucker format \begin{equation} 
\label{eq:tuck_matrix}
\widetilde{\matb{A}} =  \sum_{r_d=1}^{R^A_{d}} \dots \sum_{r_1=1}^{R^A_{1}}[\widetilde{\core{A}}]_{r_1,\dots,r_d} \widetilde{\factor{A}}_{(d,r_d)} \otimes\dots\otimes \widetilde{\factor{A}}_{(1,i_1)},
\end{equation} with  $\widetilde{\core{A}}\in\mathbb{R}^{R^A_{1}\times\dots\times R^A_{d}}$  and $\widetilde{\factor{A}}_{(i,r_i)}\in\mathbb{R}^{n_i\times R^A_{i}}$,
 and the right hand side with a vector in Tucker format 
$\widetilde{\vect{f}}$. 
In the spirit of \cite{Mantzaflaris2017}, to compute the above approximations, $f$ and the geometry coefficients are approximated by the sum of separable functions.

In this work,  this step is performed using Chebyshev polynomials.
More precisely, we use the Chebfun toolbox \cite{Driscoll2014}, which is suited for $d=3$, and in particular the \textsf{chebfun3f} function \cite{Dolgov2021} that computes low-rank approximations of trivariate functions. 
We report some details in \ref{sec:appendix2}, referring to the original papers for an exhaustive description.

In conclusion, our problem is  to find a  vector   $\widetilde{\vect{x}}$ in Tucker format that (approximately) solves
 \begin{equation}
 \label{eq:approx_ls}
 \widetilde{\matb{A}}\widetilde{\vect{x}}=\widetilde{\vect{f}}.
 \end{equation}
We remark that, since the factor matrices $\widetilde{\matb{A}}_{(i,r_i)}$, $i = 1,\ldots,d$, appearing in \eqref{eq:tuck_matrix} are banded with bandwidth $p$, then the memory required to store  $\widetilde{\matb{A}}$  is $O(d p n R_A + R_A^d)$, where, here and throughout, $R_A:=\max_{k=1,\dots,d}R^A_{k}$.
 Moreover, thanks to \eqref{eq:tt_product}, the computational cost to multiply $\widetilde{\matb{A}}$ by a vector in Tucker format is $O(d p n r R_A + r^d R_A^d )$ FLOPs, where $r$ denotes the maximum of the multilinear rank of the considered vector.

  \begin{rmk}
Other kind of boundary conditions  can be also  handled in the
low-rank setting above. For example,
if a non-homogeneous Dirichlet boundary condition $u=g$ is imposed on
$\partial{\Omega}$,  the linear system to be solved reads
\begin{equation}
  \label{eq:nonhomog-dirichlet-system}
  \widetilde{\matb{A}}\widetilde{\vect{x}}=\widetilde{\vect{f}}-
  \widetilde{\matb{A}}_{\partial \Omega}\widetilde{\vect{g}},
\end{equation}
where  $\widetilde{\matb{A}}$, as in \eqref{eq:tuck_matrix}--\eqref{eq:approx_ls}, is the matrix representing the
bilinear form 
$a(\cdot, \cdot)    $ on the basis functions that vanish on
$\partial{\Omega}$,
$\widetilde{\matb{A}}_{\partial \Omega}$  represents
$a(\cdot, \cdot)    $ on the trial basis functions that vanish on
$\partial{\Omega}$ and test basis functions whose support intersect
$\partial{\Omega}$ , $\widetilde{\vect{g}}$ is the vector of
degrees-of-freedom of (an approximation of) $g$, and  we use a
low-rank approach for the right-hand side of \eqref{eq:nonhomog-dirichlet-system}.
 \end{rmk}

\section{Low-rank linear solver}
\label{sec:low_rank_solver}
 In order to find a solution $\widetilde{\vect{x}}$ of \eqref{eq:approx_ls} in Tucker format, we  present   in Section \ref{sec:tpcg} a suited iterative solver in which each iterate is a vector in Tucker format. A fundamental role is played by the truncation operators, reviewed in Section \ref{sec:truncation}, that reduces the multilinear rank of a vector in Tucker format. Another important ingredient, that we introduce in Section \ref{sec:preconditioner} is  a novel  preconditioning strategy compatible with the Tucker format. 
 
As we will see in the following, the computational cost of each iteration of our solver is (almost) linear in the number of univariate dofs.
   
\subsection{Truncation operators}
\label{sec:truncation}
 
The sum of two Tucker vectors and the multiplication of a Tucker matrix  by a Tucker vector are operations that increase the multilinear rank (see Section \ref{sec:tensor_calculus}). For this reason it is fundamental to have an efficient strategy that allows to compress a vector in Tucker format to another vector in Tucker format with lower multilinear rank.
 In this section we review two related algorithms that perform this task: the first one is the relative tolerance truncation strategy and it is based on \cite{Oseledets2009} while the second one is the dynamic truncation strategy  from \cite{Zander2013,Matthies2012}. 
 
 Note that a straightforward application of the ST-HOSVD algorithm to the full tensor associated to a Tucker vector of dimension ${n_1\dots n_d}=N_{dof}$ is unfeasible in most of the interesting cases as  the total cost is    $O(n^{d+1})$ FLOPs, which is too high in our context.
 
 \subsubsection{Relative tolerance truncation} 
The strategy that we present in this section is based on   Algorithm 3 of  \cite{Oseledets2009}.
Let  $\epsilon>0$ a given tolerance. The truncation operator $\trunc{{T}^{rel}}$   compresses a Tucker vector $\vect{y}\in\mathbb{R}^{n_1\dots n_d}$ 
to a Tucker vector $\vect{\widetilde{y}}=\trunc{{T}^{rel}}(\vect{y},\epsilon)\in\mathbb{R}^{n_1\dots n_d}$ such that $\vect{\widetilde{y}}$ satisfies
\begin{equation}
\label{eq:rel_estimate}
\|\vect{y}-\widetilde{\vect{y}}\|_2\leq\epsilon\|\vect{y}\|_2.
\end{equation} 
We describe here the procedure to find  $\widetilde{\vect{y}}$.
Let    $\Y:=\core{Y}\times_d\factor{Y}_d\dots\times_1\factor{Y}_1\in\mathbb{R}^{n_1\times\dots\times n_d}$ be the Tucker tensor that represents $\vect{y}$,   with    core tensor $\core{Y}\in \mathbb{R}^{r_1\times\dots\times r_d}$ and factor matrices $\factor{Y}_i\in\mathbb{R}^{n_i\times r_i}$ for $i=1,\dots,d$.
The first step of the truncation procedure is the computation of the QR factorizations of the matrices $\factor{Y}_i$ for $i=1,\dots,d$:   we find  $d$  orthogonal matrices  $\matb{Q}_i\in\mathbb{R}^{n_i\times r_i}$ and $d$ upper triangular matrices $\matb{R}_i\in\mathbb{R}^{r_i\times r_i}$ such that $\factor{Y}_i=\matb{Q}_i\matb{R}_i$. Then we form explicitly the tensor $\mathcal{Z}:=\core{Y}\times_d \factor{R}_d\dots\times_1\factor{R}_1\in\mathbb{R}^{r_1\times\dots\times r_d}$
 and we compute  the Tucker approximation $$\widetilde{\mathcal{Z}}=\widetilde{\core{Y}}\times_d\matb{S}_d\dots\times_1\matb{S}_1=\textsc{sthosvd}(\mathcal{Z},\epsilon)$$ where $\widetilde{\core{Y}}\in \mathbb{R}^{s_1\times \dots \times s_d}$ is the new core tensor and $\matb{S}_i\in\mathbb{R}^{r_i\times s_i}$ are orthogonal matrices with $s_i\leq r_i$.
 We define 
 $
 \widetilde{\factor{Y}}_i:= \matb{Q}_i\matb{S}_i\in\mathbb{R}^{n_i\times s_i} 
 $ and $\widetilde{\mathcal{Y}}:=\widetilde{\core{Y}}\times_d\widetilde{\factor{Y}}_d  \dots\times_1\widetilde{\matb{Y}}_1$.
Finally, the truncated Tucker vector  $\vect{\widetilde{y}}=\trunc{{T}}^{rel}(\vect{y},\epsilon)$ is defined as the vector associated to   $\widetilde{\mathcal{Y}}$. Note that   \eqref{eq:rel_estimate} is satisfied. 
We summarize the procedure just described in Algorithm \ref{al:relative_truncation}.  
 
 \paragraph{Computational cost}   The computational cost of each of the $d$ QR decompositions of Step 1 is $O(nr^2)$ FLOPs. Step 2 has a overall computational cost of $O(dr^{d+1})$ FLOPs \cite{Vannieuwenhoven2012}.
  The ST-HOSVD at Step 3 yields a computational cost of $O(r^{d+1})$ FLOPs. Thus, the overall computational cost of Algorithm 1 is $O(dnr^2+ r^{d+1})$ FLOPs.

 \begin{algorithm}[H]

\caption{Relative tolerance truncation}\label{al:relative_truncation}
 \hspace*{\algorithmicindent} \textbf{Input}: Tucker vector  $\vect{y}=(\factor{Y}_d \otimes\dots\otimes\factor{Y}_1)\text{vec}(\core{Y})$    and the relative tolerance $\epsilon>0$. \\
 \hspace*{\algorithmicindent} \textbf{Output}: Truncated Tucker vector $\widetilde{\vect{y}}$ such that $\|\vect{y}-\widetilde{\vect{y}}\|_2\leq\epsilon\|\vect{y}\|_2$.
 \begin{algorithmic}[1]
  \State Compute QR decompositions $\factor{Y}_i=\matb{Q}_i\matb{R}_i$ for $i=1,\dots,d$;
  \State Compute $\mathcal{Z} =\core{Y}\times_d \matb{R}_d\dots\times_1\matb{R}_1$;
  \vskip 1mm
  \State Compute $\widetilde{\mathcal{Z}}=\widetilde{\core{Y}}\times_d\matb{S}_d\dots\times_1\matb{S}_1=\textsc{sthosvd}(\mathcal{Z},\epsilon)$;
  \State Compute $\widetilde{\factor{Y}}_i= \matb{Q}_i\matb{S}_i$ for $i=1,\dots,d$; 
\end{algorithmic}
 Define the Tucker tensor $ \widetilde{\mathcal{Y}}:=\widetilde{\core{Y}}\times_d\widetilde{\factor{Y}}_d  \dots\times_1\widetilde{\factor{Y}}_1$ and $\widetilde{\vect{y}}$ as the Tucker vector associated to $\widetilde{\mathcal{Y}}$.
\end{algorithm}
 
 \subsubsection{Dynamic truncation}
The residual norm in a low-rank iterative solver with relative truncation strategy can stagnate roughly at the level of the truncation error  \cite{Kressner2011}. A possible strategy to overcome this problem, proposed  in \cite{Matthies2012,Zander2013}, is to reduce the relative tolerance, whenever stagnation occurs.
The authors  introduce  an indicator of the stagnation based on the consideration that when  the progress made by the iteration operator is prevented by truncation, then there is stagnation.
Suppose that in our iterative solver we have  $\vect{x}_{k+1}=\trunc{{T}^{rel}}(\Phi^{(k)}(\vect{x}_{k}),\epsilon)$ where   $\Phi^{(k)}$ 
 represents the process of update from $\vect{x}_k$ to $\vect{x}_{k+1}$ in absence of truncation. 
 We decompose the update of the iterate in these steps
 \begin{align*}
  \widetilde{\vect{x}}_{k+1} &=\Phi^{(k)}({\vect{x}}_{k})   &&\text{iteration without   truncation,}\\
   \Delta \widetilde{\vect{x}}_{k} & = \widetilde{\vect{x}}_{k+1}- {\vect{x}}_{k}\quad\ &&\text{proposed step, in absence of truncation,}\\
    {\vect{x}}_{k+1} & = \trunc{{T}^{rel}}(\widetilde{\vect{x}}_{k+1},\epsilon) \quad\quad &&\text{truncation of the iterate,}\\
   \Delta  {\vect{x}}_k & =  {\vect{x}}_{k+1} -  {\vect{x}}_k\quad\  &&\text{actual step, taken by the perturbed iterative process}.
   \end{align*}
 The following update ratio reflects how much of the proposed update $\Delta \vect{x}_{k} $ is effectively used 
   $$
v_k  :=\frac{ \Delta \widetilde{\vect{x}}_k \cdot \ \Delta \vect{x}_k  }{ \|\Delta \widetilde{\vect{x}}_k \|_2^2  }. 
   $$
 Fixed a threshold parameter $\delta>0$ and a reducing positive parameter $\alpha<1$, we consider the update from $\vect{x}_k$ to $\vect{x}_{k+1}$      satisfactory   if 
 \begin{equation}
 \label{eq:update_condition}
  |v_k-1|\leq\delta.
 \end{equation}
If    condition \eqref{eq:update_condition} is not satisfied,   the relative tolerance $\epsilon$ is reduced by a factor $\alpha<1$ and   the truncation of $\widetilde{\vect{x}}_{k+1}$ is performed  with  relative tolerance  $\alpha \epsilon$. The process is repeated until \eqref{eq:update_condition} is  satisfied. A minimum  relative tolerance $\epsilon_{min}$ that can be reached by this process is fixed at the beginning.

  The dynamic truncation of a vector $\vect{y}_{k}$ is denoted as $[{\vect{y}}_{k+1},\epsilon_{new}]=\trunc{{T}^{dt}}(\vect{y}_k,\widetilde{\vect{y}}_{k+1},\epsilon,\alpha, \epsilon_{\min},\delta)$, where $\widetilde{\vect{y}}_{k+1}:=\Phi^{(k)}(\vect{y}_k)$ and
 $ \epsilon_{new}$ is the last relative tolerance used for the truncation. This will be used as input starting tolerance for dynamic truncation at the next TPCG iteration (see Algorithm \ref{al:cg_tensor} in Section \ref{sec:tpcg}). 
The procedure is summarized in Algorithm \ref{al:dynamic_truncation}.

 \paragraph{Computational cost} Each iteration in the loop of Algorithm \ref{al:dynamic_truncation} has the same computational cost of Algorithm \ref{al:relative_truncation}.

  \begin{algorithm}

\caption{Dynamic truncation}\label{al:dynamic_truncation}
 \hspace*{\algorithmicindent} \textbf{Input}: {Tucker vectors ${\vect{y}}_k$ and ${\widetilde{\vect{y}}}_{k+1}:=\Phi^{(k)}({\vect{y}}_k)$, initial relative tolerance $\epsilon$, threshold $\delta$, minimum relative tolerance $\epsilon_{\min}$ and reducing factor $\alpha<1$.} \\
 \hspace*{\algorithmicindent} \textbf{Output}: Truncated Tucker vector ${\vect{y}}_{k+1}$   and $\epsilon_{new}$
 \begin{algorithmic}[1]
 \State Set ${\epsilon}_{new}=\epsilon$
  \State Compute  $\Delta \widetilde{\vect{y}}_{k}   =  \widetilde{\vect{y}}_{k+1}-{\vect{y}}_{k}$;
  \State Set $ex=0;$
  \While{$ex = 0$}
  \State Compute ${\vect{y}}_{k+1}   = \trunc{{T}^{rel}}(\widetilde{\vect{y}}_{k+1},\epsilon_{new});$ 
  \State Compute $ \Delta  {\vect{y}}_k  =  {\vect{y}}_{k+1} -  {\vect{y}}_k$;
  \State Compute   $v_k=\frac{ \Delta \widetilde{\vect{y}}_k \cdot \ \Delta \vect{y}_k  }{ \|\Delta \widetilde{\vect{y}}_k \|_2^2  }$;
  \If{$|v_k-1|<\delta$} 
  \State  Set $ex=1$;
  \ElsIf{$\alpha {\epsilon}>\epsilon_{min}$}
 \State Set ${\epsilon}_{new}=\alpha \epsilon_{new}$;
  \Else 
  \State  Set $ex=1;$
   \EndIf
     \EndWhile
\end{algorithmic} 
\end{algorithm}

\subsection{The preconditioner}
\label{sec:preconditioner}

In this section we extend the preconditioner proposed in \cite{Sangalli2016}, based on the Fast Diagonalization (FD) method \cite{Lynch1964}, to the low-rank setting. For the sake of simplicity, we give details only for the three-dimensional case, although it can be easily generalized to any dimension $d$.

The FD-based preconditioner is defined as
\begin{equation}
\label{eq:FD_prec}
\matb{P}:=\widehat{\matb{K}}_3\otimes\widehat{\matb{M}}_2\otimes\widehat{\matb{M}}_1 + \widehat{\matb{M}}_3\otimes\widehat{\matb{K}}_2\otimes\widehat{\matb{M}}_1 +\widehat{\matb{M}}_3\otimes\widehat{\matb{M}}_2\otimes\widehat{\matb{K}}_1,
\end{equation}
where $\widehat{\matb{K}}_i$ and $\widehat{\matb{M}}_i$ are the univariate stiffness and mass matrices  in the $i$-th parametric direction, respectively.
  Its  application can be performed first by computing the generalized eigendecomposition of $(\widehat{\matb{K}}_i, \widehat{\matb{M}}_i)$ for $i=1,2,3$. In this way, we find $\widehat{\matb{M}}_i$-orthogonal matrices $\matb{U}_i$ and diagonal matrices $\boldsymbol{\Lambda}_i$ for $i=1,2,3$, such that
\begin{equation}
\label{eq:factorizations}
\widehat{\matb{K}}_i\matb{U}_i=\widehat{\matb{M}}_i\matb{U}_i\boldsymbol{\Lambda}_i\quad \text{ and } \quad  \matb{U}_i^{T}\widehat{\matb{M}}_i\matb{U}_i=\mathbf{I}_{n_i},
\end{equation}
where  $\mathbf{I}_{n_i}\in\mathbb{R}^{n_i\times n_i}$ is the identity matrix.   Thus we have for $i=1,2,3$
$$
\widehat{\matb{K}}_i=\matb{U}_i^{-T}\boldsymbol{\Lambda}_i\matb{U}_i^{-1}\quad \text{ and } \quad   \widehat{\matb{M}}_i =\matb{U}_i^{-T}\matb{U}_i^{-1}
$$
and by inserting the above factorization in \eqref{eq:FD_prec}, we get
\begin{equation*}
\matb{P}= (\matb{U}_3 \otimes\matb{U}_2 \otimes \matb{U}_1)^{-T}( \mathbf{I}_{n_3}\otimes\mathbf{I}_{n_2}\otimes \boldsymbol{\Lambda}_1
 +\mathbf{I}_{n_3}\otimes\boldsymbol{\Lambda}_2\otimes \mathbf{I}_{n_1}+\boldsymbol{\Lambda}_3\otimes\mathbf{I}_{n_2}\otimes \mathbf{I}_{n_1})(\matb{U}_3 \otimes\matb{U}_2 \otimes \matb{U}_1)^{-1}.
\end{equation*}
Therefore, its inverse is
\begin{equation}
\label{eq:factorized_FD}
\matb{P}^{-1}= (\matb{U}_3 \otimes\matb{U}_2 \otimes \matb{U}_1)( \mathbf{I}_{n_3}\otimes\mathbf{I}_{n_2}\otimes \boldsymbol{\Lambda}_1
 +\mathbf{I}_{n_3}\otimes\boldsymbol{\Lambda}_2\otimes \mathbf{I}_{n_1}+\boldsymbol{\Lambda}_3\otimes\mathbf{I}_{n_2}\otimes \mathbf{I}_{n_1})^{-1}(\matb{U}_3 \otimes\matb{U}_2 \otimes \matb{U}_1)^T.
\end{equation}
 We summarize the standard FD method in Algorithm \ref{al:FD}.
\begin{algorithm}[H]
\caption{Fast Diagonalization}\label{al:FD}
\hspace*{\algorithmicindent} SETUP OF THE PRECONDITIONER:\\ 
\hspace*{\algorithmicindent} \textbf{Input}: {Pencils $(\widehat{\matb{K}}_i,\widehat{\matb{M}}_i)$ for $i=1,2,3$.  } \\
 \hspace*{\algorithmicindent} \textbf{Output}:  {Preconditioner \eqref{eq:factorized_FD}.}
\begin{algorithmic}[1] 
\State Compute the factorizations \eqref{eq:factorizations}.
\Statex 
\Statex APPLICATION OF THE PRECONDITIONER:
\Statex \textbf{Input}: {Factorizations \eqref{eq:factorizations} and right-hand side   vector $\vect{s} \in\mathbb{R}^{n_1  n_2 n_3}.$  } 
\Statex \textbf{Output} Solution of $\mathbf{P}\vect{v}=\vect{s}.$
 
\State Compute $\widetilde{\vect{s}} = (\matb{U}_3 \otimes\matb{U}_2 \otimes \matb{U}_1)^T\vect{s}$;
\State Compute $\widetilde{\vect{q}} = \left( \mathbf{I}_{n_3}\otimes\mathbf{I}_{n_2}\otimes  \boldsymbol{\Lambda}_1
 +\mathbf{I}_{n_3}\otimes \boldsymbol{\Lambda}_2\otimes \mathbf{I}_{n_1}+ \boldsymbol{\Lambda}_3\otimes\mathbf{I}_{n_2}\otimes \mathbf{I}_{n_1}\right)^{-1} \widetilde{\vect{s}};$
\State Compute $\vect{v} = (\matb{U}_3 \otimes\matb{U}_2 \otimes \matb{U}_1)\ \widetilde{\vect{q}}.$
\end{algorithmic}
\end{algorithm}
It is clear that the application of $\matb{P}$ as described in  Algorithm \ref{al:FD} is not suited to our context. 
Indeed, the inverse diagonal matrix appearing in Step 3 is not in Tucker format. Moreover, the setup and application of the preconditioner yield a computational cost that greatly exceeds the ideal $O(n)$ cost. 
 
In order to reduce the computational complexity  of both  setup and  application of the preconditioner, we  propose a technique to approximate the eigendecompositions that is not only cheaper than the exact factorization, but also reduces the cost to compute matrix-vector products with the eigenvector matrices. 
We report this strategy in \ref{sec:appendix}; a detailed analysis is postponed to a forthcoming publication. Nevertheless, we emphasize that the numerical experiments presented in Section \ref{sec:numerics} indicate the effectiveness of this approach. 

Let $\widetilde{\matb{U}}_1$, $\widetilde{\matb{U}}_2$, $\widetilde{\matb{U}}_3$  and $\widetilde{\boldsymbol{\Lambda}}_1$, $\widetilde{\boldsymbol{\Lambda}}_2$, $\widetilde{\boldsymbol{\Lambda}}_3$ represent the approximated eigenvector and eigenvalue matrices, respectively. 
The corresponding preconditioner takes the form 
\begin{align}
\matb{P}^{-1}\approx    (\widetilde{\matb{U}}_3 \otimes\widetilde{\matb{U}}_2 \otimes \widetilde{\matb{U}}_1)  \matb{D}(\widetilde{{\matb{U}}}_3 \otimes\widetilde{\matb{U}}_2 \otimes \widetilde{\matb{U}}_1)^T,
\end{align}
 where $\matb{D}:=  (\mathbf{I}_{n_3}\otimes\mathbf{I}_{n_2}\otimes \widetilde{\boldsymbol{\Lambda}}_1
 +\mathbf{I}_{n_3}\otimes\widetilde{\boldsymbol{\Lambda}}_2\otimes \mathbf{I}_{n_1}+\widetilde{\boldsymbol{\Lambda}}_3\otimes\mathbf{I}_{n_2}\otimes \mathbf{I}_{n_1})^{-1}$. 

 We build an approximation $\widetilde{\matb{D}} $ of $\matb{D}$ that can be written in Tucker format and   that satisfies
\begin{equation} \label{eq:prec_bound}
\| \matb{D}^{-1} \|_{\infty} \|\matb{D}-\widetilde{\matb{D}}\|_{\infty}  \leq \epsilon_{rel}^{prec}
\end{equation}
for a given relative tolerance $\epsilon_{rel}^{prec}>0$. Note that since
$$   \| \mathbf{I}_{N_{dof}} - \matb{D}^{-1} \widetilde{\matb{D}}\|_{\infty} \leq \| \matb{D}^{-1} \|_{\infty} \|\matb{D}-\widetilde{\matb{D}}\|_{\infty} , $$
 inequality \eqref{eq:prec_bound} guarantees in particular that all eigenvalues of $\matb{D}^{-1}\widetilde{\matb{D}}$ belong to the interval $\left[ 1 - \epsilon_{rel}^{prec}, 1 + \epsilon_{rel}^{prec} \right] $.
 
Our approach is based on the approximation of the function $g(\lambda) = 1/\lambda$ using a linear combination of exponential functions. We first recall the following result \cite{Braess2005,Kressner2010}.

\begin{prop} 
\label{lemma:approx_diag}
 Let $s_R(\lambda):=\sum_{j=1}^R  \omega_j\exp(-\alpha_j \lambda)$, with $\lambda, \alpha_j,\omega_j\in\mathbb{R}$ and $R\in\mathbb{N}$. Then, there is a choice of $\alpha_j>0$ and $\omega_j>0$ for $j=1,\ldots,R$, that depends on $R$ and $M$, such that
\begin{align}
\label{eq:estimate_approx_diag}
 E_R(M):=\sup_{\lambda\in[1,M]}\left|\frac{1}{\lambda}-s_R(\lambda) \right|\leq 16 \exp\left( \frac{-R\pi^2}{\log(8M)}\right).
\end{align}
 \end{prop}
The explicit values of $\alpha_j$ and $\omega_j$ of  Proposition \ref{lemma:approx_diag} are not known. 
However, good approximations of these parameters and of $E_R(M)$ are provided in \cite{Hackbusch2019} for a wide range of values of $R$ and $M$.
 We define $\lambda_{\min}:=\min (\widetilde{\boldsymbol{\Lambda}}_1) +\min(\widetilde{\boldsymbol{\Lambda}}_2)+\min(\widetilde{\boldsymbol{\Lambda}}_3)$ and  $\lambda_{\max}:=\max (\widetilde{\boldsymbol{\Lambda}}_1) +\max(\widetilde{\boldsymbol{\Lambda}}_2)+\max(\widetilde{\boldsymbol{\Lambda}}_3)$. Given a relative tolerance $\epsilon_{rel}^{prec}$, we define  
\begin{equation}
\label{eq:estimate_prec_rel}
  M_P:=\frac{\lambda_{\max}}{\lambda_{\min}}\quad\text{   and   }  \quad R_P:=\min\left\{R \ \big| \ E_R(M_P)\leq \frac{\epsilon_{rel}^{prec}}{M_P}\right\}.
 \end{equation}
We take $M=M_P$ and $R=R_P$ and we consider the associated   (approximated) optimal parameters $\omega_j$ and $\alpha_j$ for $j=1,\ldots,R_P$ from Proposition \ref{lemma:approx_diag}. Then we define
\begin{equation}
\label{eq:approx_diagonal}
\widetilde{ \matb{D}} :=\frac{1}{\lambda_{\min}}\sum^{R_P}_{j=1}\omega_j \matb{D}_{(3,j)}\otimes \matb{D}_{(2,j)}\otimes\matb{D}_{(1,j)}
 \end{equation}
  with $\matb{D}_{(i,j)}$ diagonal matrices with diagonal entries equal to
\begin{align}
\label{eq:diag_factor}
[\matb{D}_{(i,j)}]_{m,m}&:=\exp{\left(-\frac{\alpha_j}{\lambda_{\min}}[\widetilde{\boldsymbol{\Lambda}}_i]_{m,m}\right)}\quad \text{for }j=1,\dots R_P,\ m=1,\dots,n_i \text{ and } i=1,2,3.
\end{align}
Thanks to \eqref{eq:estimate_approx_diag} and the definition of $R_P$, we have
$$ \|\matb{D}-\widetilde{\matb{D}} \|_{\infty} \leq  \sup_{\lambda\in[\lambda_{\min},\lambda_{\max}]}\left|\frac{1}{\lambda}-\frac{1}{\lambda_{\min}} s_{R_P}(\lambda) \right| =  \displaystyle \frac{1}{\lambda_{\min}} \sup_{\lambda\in[1,M_P]}\left|\frac{1}{\lambda}-s_{R_P}(\lambda) \right|  = \frac{1}{\lambda_{\min}} E_{R_P}(M_P).$$
Since $\|\matb{D}^{-1}\|_{\infty} = \lambda_{\max} $, it holds 
$$ \|\matb{D}^{-1}\|_{\infty} \|\matb{D}-\widetilde{\matb{D}} \|_{\infty} \leq M_P E_{R_P}(M_P) \leq \epsilon_{rel}^{prec}. $$
Finally, our preconditioner is defined as
\begin{equation}
\label{eq:prec}
\widetilde{\matb{P}}^{-1}:=  \frac{1}{\lambda_{\min}} (\widetilde{\matb{U}}_3 \otimes\widetilde{\matb{U}}_2 \otimes \widetilde{\matb{U}}_1) \left( \sum^{R_P}_{j=1}\omega_j  \matb{D}_{(3,j)}\otimes \matb{D}_{(2,j)} \otimes \matb{D}_{(1,j)}\right) (\widetilde{\matb{U}}_3 \otimes\widetilde{\matb{U}}_2 \otimes \widetilde{\matb{U}}_1)^T,
\end{equation}
which can be written as a Tucker matrix as 
$$
\widetilde{\matb{P}}^{-1}:=  \sum^{R_P}_{j_3=1} \sum^{R_P}_{j_2=1} \sum^{R_P}_{j_1=1} [\core{P}]_{j_1,j_2,j_3} \matb{P}_{(3,j_3)}\otimes \matb{P}_{(2,j_2)} \otimes \matb{P}_{(1,j_1)},
$$
where $\core{P}\in\mathbb{R}^{{R_P}\times {R_P} \times {R_P}} $ is the diagonal tensor whose diagonal entries are 
\begin{equation}
\label{eq:core_tens_prec}
[\core{P} ]_{j,j,j}= \frac{\omega_{j}}{\lambda_{\min}}, \qquad  j = 1,\ldots, R_P,
\end{equation}
and where $\matb{P}_{(i,j)}:=  \widetilde{\matb{U}}_i \matb{D}_{(i,j)} \widetilde{\matb{U}}_{i}^T$ for $i=1,2,3$ and $j=1,\dots,R_P$.
%
%
The setup and application of our preconditioner are summarized in Algorithm \ref{al:approx_FD}.

\begin{algorithm}
\caption{Low-rank Fast Diagonalization}\label{al:approx_FD}
 \hspace*{\algorithmicindent} SETUP OF THE PRECONDITIONER:\\
 \hspace*{\algorithmicindent}\textbf{Input}: {$\epsilon_{rel}^{prec}$ relative tolerance of the preconditioner.} \\
  \hspace*{\algorithmicindent}\textbf{Output}: {Preconditioner $\widetilde{\matb{P}}$ in Tucker format as \eqref{eq:prec}.}
\begin{algorithmic}[1]
\State 
Find  $\widetilde{\matb{U}}_1$,$\widetilde{\matb{U}}_2$ and $\widetilde{\matb{U}}_3$  as described in  \ref{sec:appendix};
\State Define $M_P:=\lambda_{\max}/\lambda_{\min}$ and find $R_P$ such that $E_R(M_P)\leq \frac{\epsilon_{rel}^{prec}}{M_P}$;
\State Compute   the factors $\matb{D}_{(i,j)}$ for $j=1,\dots,R_P$ and $i=1,2,3$ as in \eqref{eq:diag_factor};
\State Define $\core{P}$ as in \eqref{eq:core_tens_prec}.
\Statex 
\Statex APPLICATION OF THE PRECONDITIONER:
\Statex \textbf{Input}: {Tucker  vector $\vect{s}$ represented by $\mathcal{S}=\core{S}\times_3\factor{S}_3\times_2\factor{S}_2\times_1\factor{S}_1$ where  $\core{S}\in\mathbb{R}^{r_1\times r_2\times r_3}$,  $R_P$ rank of the approximation \eqref{eq:approx_diagonal}}.
 \Statex  \textbf{Output}:  {Solution of $\widetilde{\matb{P}}\vect{v}=\vect{s}$ in Tucker format.}
\State  Compute $\core{V}=\core{P}\kron\core{S}$;
\State Compute $\matb{Z}_i = \widetilde{\matb{U}}_i^T\factor{S}_i$ for  $i=1,2,3$;
\State Compute $\matb{W}_{(i,j)} =\matb{D}_{(i,j)}\matb{Z}_i$ for $j=1,\dots,R_P$ and $i=1,2,3$; 
\State Compute $\matb{Y}_{(i,j)} = \widetilde{\matb{U}}_i\matb{W}_{(i,j)}$ for $j=1,\dots,R_P$ and $i=1,2,3$.     
\State Define $\vect{v}$ as the Tucker vector  represented by $\mathcal{V}:=\core{V}\times_3 \factor{V}_3 \times_2 \factor{V}_2\times_1\factor{V}_1$ where $\factor{V}_i:=[\matb{Y}_{(i,1)}, \dots, \matb{Y}_{(i,R_P)}]$ for $i=1,2,3.$\end{algorithmic}
\end{algorithm}
\paragraph{Computational cost}  
The computational cost of the setup of the preconditioner is negligible, see  \ref{sec:appendix}. 
As for the application cost, assume to apply $\widetilde{\matb{P}}^{-1}$ to a vector in Tucker format having maximum of the multilinear rank $r$. Then Step 5 has a  cost of $O(r^3 R_P^3)$ FLOPs,   while each matrix-matrix multiplication of Step 6 and Step 8 requires $O(n r (\log(n)+p))$ FLOPs, thanks to the FFT (see  \ref{sec:appendix}). Finally, Step 7, being just a diagonal scaling, has a negligible cost. 
Generalizing to dimension $d$, the total cost for Algorithm \ref{al:approx_FD} is $O( d n r R_P (\log(n)+p)+ r^d R_P^d )$ FLOPs.

 \subsection{Truncated PCG method}
 \label{sec:tpcg}
The low-rank iterative solver that we choose for solving system \eqref{eq:approx_ls} is the   Truncated Preconditioned Conjugate Gradient (TPCG) method \cite{Kressner2011}, that we report in Algorithm \ref{al:cg_tensor}.
All the vectors and matrices present in the algorithm have to be intended as Tucker vectors and Tucker matrices, respectively, and thus   the operations of matrix-vector product and scalar product exploit this assumption for an efficient computation (see Section \ref{sec:tensor_calculus}). 
The application of our preconditioner to a Tucker vector  increases its multilinear rank, thus compared to \cite{Kressner2011} we apply one additional truncation step, and  precisely at Step 13 of Algorithm \ref{al:cg_tensor}.  
Note also that we compute the residual directly as $\vect{\widetilde{f}} - \matb{\widetilde{A}}\vect{x}_{k}$ because, as observed in \cite{Kressner2011}, the recursive computation of the residual can lead to stagnation.  
 We apply the dynamic truncation to the  iterate $\vect{x}_k$ and relative truncation to the other iteration vectors $\vect{r}_k, \vect{z}_k , \vect{p}_k$ and $\vect{q}_k$. The tolerance of latter truncations should be chosen carefully, in order to keep the multilinear ranks of  $\vect{r}_k, \vect{z}_k , \vect{p}_k$ comparable to the multilinear rank of $\vect{x}_k$. 
Indeed, in our numerical experience if we set the relative truncation tolerance of those iterates equal to the one used for the truncation of $\vect{x}_k$, then the multilinear rank of those iterates becomes significantly higher than one of  $\vect{x}_k$. Thus, inspired by the numerical experiments in \cite{Palitta2021} and by \cite{Simoncini2003}, we choose a relative tolerance for the truncation of $\vect{r}_k, \vect{z}_k , \vect{p}_k$ and $\vect{q}_k$ equal to 
\begin{equation}
\label{eq:rel_tol_iterates}
\eta_k:=\beta\frac{tol\|\widetilde{\vect{f}}\|_2}{\|\vect{r}_k\|_2}
\end{equation}
where $0<\beta<1$ and $tol$ is the tolerance of  TPCG. At the beginning of the iterative process the residual is large and thus the truncation of the iterates is more aggressive, while when convergence is approached, the residual is small and the truncation is relaxed, resulting in a limited multilinear rank.

\begin{algorithm}

\caption{TPCG}\label{al:cg_tensor}
  \hspace*{\algorithmicindent} \textbf{Input}: { Linear system matrix   $\matb{\widetilde{A}}$  and preconditioner $\widetilde{\matb{P}}$ in Tucker format,   right-hand side Tucker vector  $\vect{\widetilde{f}}$,  initial guess $\vect{x}_0$ in Tucker format,   TPCG tolerance $tol>0$, parameter $\beta$ for the relative truncation,  
parameters for the dynamic truncation: starting relative tolerance $\epsilon_0$, reducing factor $\alpha$, minimum relative tolerance $\epsilon_{\min}$, threshold $\delta$.} \\
 \hspace*{\algorithmicindent} \textbf{Output}: Low-rank solution $\widetilde{\vect{x}}$  of $\matb{\widetilde{A}} \widetilde{\vect{x}} =\vect{\widetilde{f}}.$
 \begin{algorithmic}[1]
  \State $\vect{r}_0=\vect{\widetilde{f}}-\matb{\widetilde{A}}\vect{x}_0$;
  \State $\eta_0 = \beta\frac{tol\|\vect{f}\|_2}{\|\vect{r}_0\|_2}$;
  \State $\vect{z}_0=\trunc{{T}^{rel}}(\widetilde{\matb{P}}^{-1}\vect{r}_0,\eta_0)$;
  \State $\vect{p}_0=\vect{z}_0$;
  \State $\vect{q}_0=\trunc{{T}^{rel}}(\matb{\widetilde{A}}\vect{p}_0,\eta_0)$;
  \State ${\xi}_0=\vect{p}_0\cdot\vect{q}_0$;
  \State $k=0$
 \While{$\|\vect{r}_{k}\|_2>tol$}{\\
        \qquad  $\omega_k = \frac{\vect{r}_k \cdot \vect{p}_k }{{\xi}_k}$;\\ 
         \qquad  $[\vect{x}_{k+1},\epsilon_{k+1}]={    \trunc{{T}^{dt}}}(\vect{x}_k, \vect{x}_k+\omega_k\vect{p}_k,\epsilon_k, \alpha, \epsilon_{\min},\delta)$; \\  
         \qquad  $\vect{r}_{k+1}={  \trunc{{T}^{rel}}}(\vect{\widetilde{f}} - \matb{\widetilde{A}}\vect{x}_{k+1},\eta_{k})$ ; \\
         \qquad $\eta_{k+1}=\beta\frac{tol\|\vect{f}\|_2}{\|\vect{r}_{k+1}\|_2}$;\\
         \qquad  $\vect{z}_{k+1}=   \trunc{{T}^{rel}}(\widetilde{\matb{P}}^{-1}  \vect{r}_{k+1},\eta_{k+1})$; \\
         \qquad  $\beta_k=-\frac{\vect{z}_{k+1}\cdot\vect{q}_k} {\vect{ \xi}_{k}}$;\\
        \qquad    $\vect{p}_{k+1}={    \trunc{{T}^{rel}}}(\vect{z_{k+1}}+\beta_k\vect{p}_k, \eta_{k+1})$;  \\
        \qquad  $\vect{q}_{k+1}={   \trunc{{T}^{rel}}}(\matb{\widetilde{A}}\vect{p}_{k+1},\eta_{k+1})$;  \\
        \qquad ${\xi}_{k+1}=\vect{p}_{k+1}\cdot \vect{q}_{k+1}$\\
        \qquad   $k=k+1$;
  }\EndWhile
\State $\widetilde{\vect{x}}=\vect{x}_{k}.$
\end{algorithmic}
\end{algorithm}

\paragraph{Computational cost}   
Here we summarize the computational cost of each iteration of Algorithm \ref{al:cg_tensor}. In this regard, the main efforts are represented by the residual computation (Step 11), the application of the preconditioner (Step 13), and the computation of $\matb{\widetilde{A}}\vect{p}_{k+1}$ (Step 16), as well as the corresponding truncation steps. Let $r_{\vect{x}}$, $r_{\vect{r}}$ and $r_{\vect{p}}$ denote respectively the maximums of the multilinear ranks of $\vect{x}_{k+1}$, $\vect{r}_{k+1}$ and $\vect{p}_{k+1}$.
Then the maximums of the multilinear ranks of $\matb{\widetilde{A}}\vect{x}_{k+1}$, $\matb{\widetilde{P}}^{-1}\vect{r}_{k+1}$ and $\matb{\widetilde{A}}\vect{p}_{k+1}$ are bounded  respectively by $R_A r_{\vect{x}}$, $R_P r_{\vect{r}}$ and $R_A r_{\vect{p}}$.
Then, according to the cost analysis reported in sections \ref{sec:model_problem}, \ref{sec:truncation} and \ref{sec:preconditioner}, the number of FLOPs required by these steps can be bounded as follows.

\begin{table}[H]
\begin{tabular}{llll}
Residual computation: & $O(d p n r_{\vect{x}} R_A  + r_{\vect{x}}^d R_A^d )$ & Truncation: & $O( d n r_{\vect{x}}^2 R_A^2 +  r_{\vect{x}}^{d+1} R_A^{d+1} )$ \\
Preconditioner application: & $O(d n r_{\vect{r}} R_P (\log(n) + p) + R_P^d r_{\vect{r}}^d)$ & Truncation: & $O(d n r_{\vect{r}}^2 R_P^2 +  r_{\vect{r}}^{d+1} R_P^{d+1})$ \\
Computation of $\matb{\widetilde{A}}\vect{p}_{k+1}$: & $O(d p n r_{\vect{p}} R_A + r_{\vect{p}}^d R_A^d )$ & 
Truncation: & $O(d n r_{\vect{p}}^2 R_A^2 +  r_{\vect{p}}^{d+1} R_A^{d+1})$
\end{tabular}
\end{table}

Here we assumed that the maximum of the multilinear rank of $\vect{\widetilde{f}} - \matb{\widetilde{A}}\vect{x}_{k+1}$ is comparable to that of $\matb{\widetilde{A}}\vect{x}_{k+1}$ (which is reasonable since, in order to have $\matb{\widetilde{A}}\vect{x}_{k+1} \approx \vect{\widetilde{f}}$ the rank of $\matb{\widetilde{A}}\vect{x}_{k+1}$ is likely comparable or higher than the rank of $\vect{\widetilde{f}}$).  
The above analysis demonstrates that  the  computational cost is (almost) linear with respect to $n$. Moreover, we can see that the main effort is likely represented by the truncation steps. This is in agreement with our numerical experience.

\section{Numerical experiments}
\label{sec:numerics}
In this section we propose some numerical experiments to assess the performance of our low-rank strategy.
 The tests are performed using   Matlab R2021a  on an 16-Core Intel Xeon W, running at 3.20 GHz, and with  384 GB of RAM. Isogeometric   discretizations are handled with   the GeoPDEs toolbox \cite{Vazquez2016}, while operations involving tensors are managed with the Tensorlab toolbox \cite{Sorber2014}. For the approximation of the linear system matrix and right-hand side in Tucker format, we make use of the Chebfun toolbox \cite{Driscoll2014}.
 
We apply the low-rank solver described in the previous sections to the Poisson problem on two computational domains. Moreover, in Section \ref{sec:le}, we adapt this strategy to the compressible linear elasticity problem and we provide a numerical experiment to assess its good behavior. In all cases, we consider a dyadic refinement of the domain yielding  $n_{el}:=2^{l}$ elements in each parametric direction, for different values of the discretization level $l$.

The TPCG method presented in Algorithm \ref{al:cg_tensor} is employed as linear iterative solver.  The relative tolerance for the approximation of the system matrix and the right-hand side  (see \ref{sec:appendix2})  is set equal to $\epsilon=\max(tol\cdot 10^{-1},10^{-12})$, where $tol$ is the relative tolerance of the TPCG method, whose choice depends on the considered problem. In all tests we fix  the initial guess as the zero vector $\vect{x}_0=\vect{0}$. For the dynamic truncation (Algorithm \ref{al:dynamic_truncation}), we fix the initial relative tolerance $\epsilon_0$ equal to $10^{-1}$ and the minimum relative tolerance   $\epsilon_{\min}$ to $tol\| \widetilde{\vect{f}} \|_210^{-1}$, where $\widetilde{\vect{f}}$ is the right-hand side vector of \eqref{eq:approx_ls}. The reducing factor $\alpha$ is set equal to $0.5$, while the threshold $\delta$ is chosen to be  $10^{-3}$. The relative tolerance parameter $\beta$ in \eqref{eq:rel_tol_iterates} is fixed equal to $10^{-1}$. These choices of parameters yield good performances in the convergence of TPCG method. 
  
Let $\widetilde{\vect{x}}$ be the Tucker solution of the Poisson problem obtained
with the TPCG algorithm, and assume it is represented by the Tucker tensor
$\widetilde{\mathcal{X}}=\core{X}\times_3\factor{X}_3\times_2\factor{X}_2\times_1\factor{X}_1$
with $\core{X}\in\mathbb{R}^{r_1\times r_2\times r_3}$ and
$\factor{X}_i\in\mathbb{R}^{n_i\times r_i}$ for $i=1,2,3$. We define
the memory compression percentage of the low-rank solution with
respect to the full solution as 
\begin{equation}
\label{eq:mem_comp}
\text{memory compression} = \frac{r_1r_2r_3 +
  r_1n_1+r_2n_2+r_3n_3}{n_1 n_2 n_3}\cdot 100. 
\end{equation}
Note that here we are not taking into account the memory compression for the stiffness matrix (which was already analyzed and tested extensively in \cite{Mantzaflaris2017,Hofreither2018,Bunger2020}) and of the right-hand side.

In the numerical experiments we report the maximum of the multilinear rank of the solution. We remark that, thanks to the relaxation strategy on the truncation tolerance of the iterates $\vect{r}_k, \vect{z}_k , \vect{p}_k$ and $\vect{q}_k$, when convergence approaches (see Section \ref{sec:truncation}) the {maxima } of their multilinear ranks is similar to the one of $\vect{x}_k$, and therefore their values are not reported.

In all tests related to the Poisson problem we use \eqref{eq:prec} as a preconditioner with $\epsilon_{rel}^{prec}=10^{-1}$. The values of $M_P$ and the corresponding  ranks $R_P$, as defined in \eqref{eq:estimate_prec_rel},   that we found for different degrees $p$ and number of elements $n_{el}$ for Poisson problem are reported in  Table \ref{tab:ranks_table}. We emphasize that the values of $R_P$ are relatively small and show only a mild dependence with respect to the number of elements and the degree. 

We recall that $N_{dof}$ represents the total number of degrees of freedom, that is the dimension of the B-spline space in \eqref{eq:spline_space}.
 
 {\renewcommand\arraystretch{1.4} 
\begin{table}
\begin{center}
\begin{tabular}{|c|c|c|c|c|}
\hline
 & \multicolumn{4}{|c|}{ \ Preconditioner $\widetilde{\matb{P}}$: values of   $M_P$ /  $R_P$ } \\
 \hline
$n_{el}$ & $p=2$  & $p=3$ &$p=4$  &$p=5$   \\
\hline 
128 &  $1.6\cdot 10^4$ / 11 & $2.3\cdot 10^4$ / 12  &   $4.0\cdot 10^4$ / 13 & $6.5\cdot 10^4$ / 13\\
\hline
256 &  $6.6\cdot 10^4$   / 13  & $9.5\cdot 10^4$ / 14  &   $1.6\cdot 10^5$ / 15 &  $2.6\cdot 10^5$/ 16\\
\hline
512 &    $2.6\cdot 10^5$  / 16  & $3.8\cdot 10^5$ / 17  & $6.4\cdot 10^5$ / 18  &  $1.0\cdot 10^6$ / 19 \\ 
\hline
1024 &  $1.0\cdot 10^6$ / 19      & $1.5\cdot 10^6$ / 19  &  $2.6\cdot 10^6$ / 21 & $4.1\cdot 10^6$  / 22 \\ 
\hline
\end{tabular}
\caption{Values of $M_P:=\frac{\lambda_{\max}}{\lambda_{\min}}$ and corresponding preconditioner rank $R_P$ with $\epsilon_{rel}^{prec}=10^{-1}$ for different number of elements $n_{el}$ and degrees $p$ for the Poisson problem.}
\label{tab:ranks_table}
\end{center}
\end{table}}

\subsection{Thick quarter of annulus domain}
\label{sec:tq_test}
\begin{figure}[H]
\begin{center}
 \includegraphics[scale=0.45]{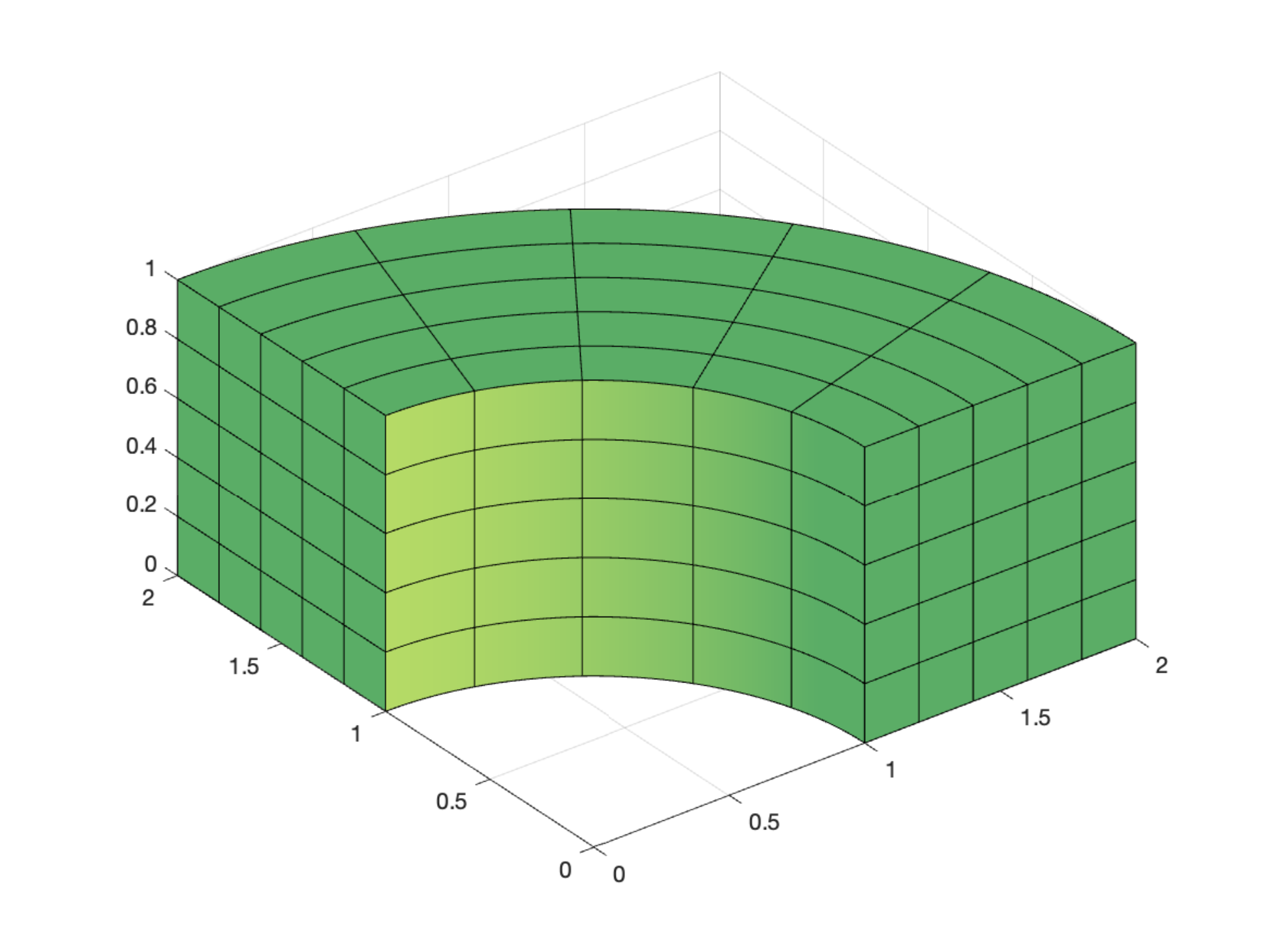}
 \caption{Thick quarter of annulus domain}
 \label{fig:thick-quarter}
 \end{center}
 \end{figure}

In this test we consider as computational domain a thick quarter of annulus, represented in  Figure \ref{fig:thick-quarter}.

We consider the manufactured solution $u(x,y,z)= (x^2+y^2-1)(x^2+y^2-4)\sin (\pi z)\sin(7xy)$ with homogeneous Dirichlet boundary condition.  
While strictly speaking this function is not low-rank, it can be effectively approximated with a low-rank solution, as we will see in the following.

The multilinear ranks of the coefficients of the  stiffness matrix in equation \eqref{eq:matrix_CP} that we found using the technique described in  \ref{sec:appendix2} are
$$
\begin{bmatrix}
(R_1^{(1,1)},R_2^{(1,1)},R_3^{(1,1)}) &  (R_1^{(1,2)},R_2^{(1,2)},R_3^{(1,2)}) &  (R_1^{(1,3)},R_2^{(1,3)},R_3^{(1,3)}) \\
 (R_1^{(2,1)},R_2^{(2,1)},R_3^{(2,1)}) & (R_1^{(2,2)},R_2^{(2,2)},R_3^{(2,2)}) & ( R_1^{(2,3)},R_2^{(2,3)
},R_3^{(2,3)})\\
(R_1^{(3,1)},R_2^{(3,1)},R_3^{(3,1)}) & (R_1^{(3,2)},R_2^{(3,2)},R_3^{(3,2)}) & (R_1^{(3,3)},R_2^{(3,3)},R_3^{(3,3)})
\end{bmatrix}=\begin{bmatrix}
(1,1,1) &  \vect{0} &  \vect{0} \\
 \vect{0} & (1,1,1) & \vect{0}\\
  \vect{0} &  \vect{0} & (1,1,1)
\end{bmatrix},
$$
where $ \vect{0}=(0,0,0)$. Thus, referring to equation \eqref{eq:tuck_matrix}, we have
$$
  \left(R_1^{A}, R_2^A, R_3^A\right) = (3,3,3),
$$
indicating that the domain has a natural tensor product structure.

\begin{figure}[H]
\centering
\subfloat[][$L^2$ error in the thick quarter domain.\label{fig:L2_quarter}]
{\includegraphics[scale=0.45]{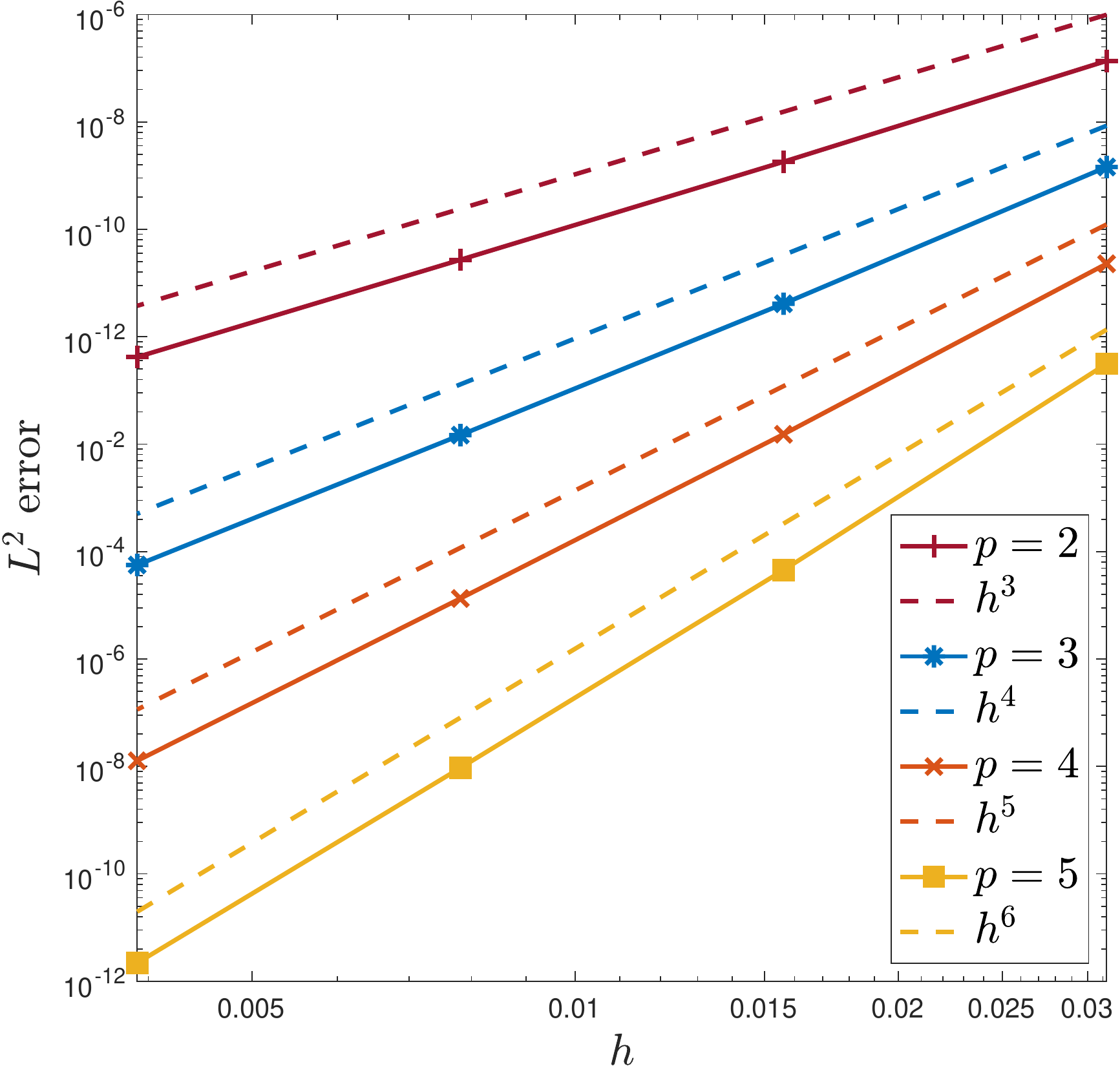}}\quad
\subfloat[][$H^1$ error in the thick quarter domain.\label{fig:H1_quarter}]
{\includegraphics[scale=0.45]{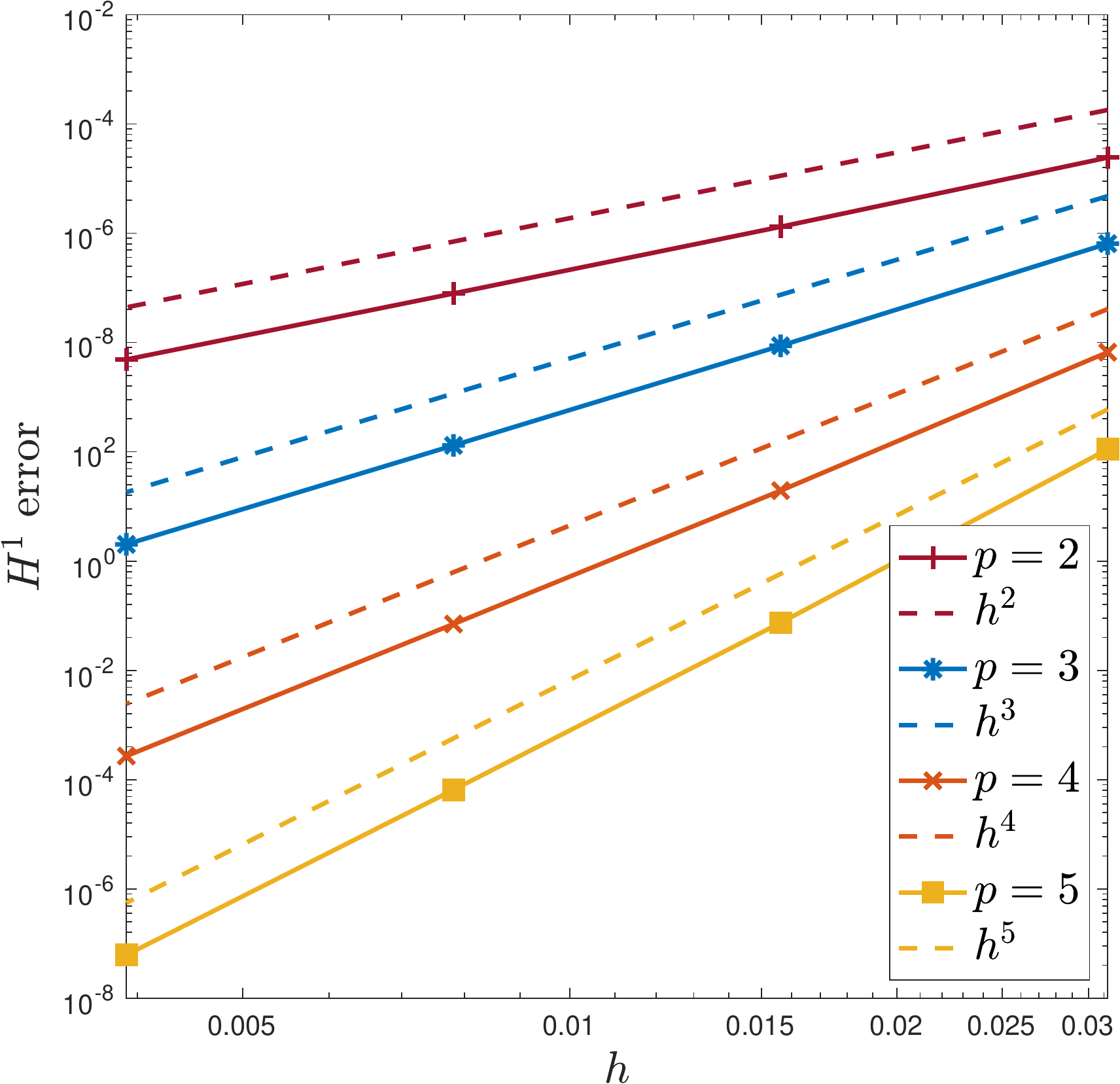}}
 \caption{Errors in the thick quarter domain with variable TPCG tolerance.}
 \label{fig:quarter_error}
\end{figure}

\begin{figure}[H]
 \centering  
  \subfloat[][Maximum of the multilinear rank of the solution.\label{fig:rank_TQ_error}]  
    {\includegraphics[scale=0.45]{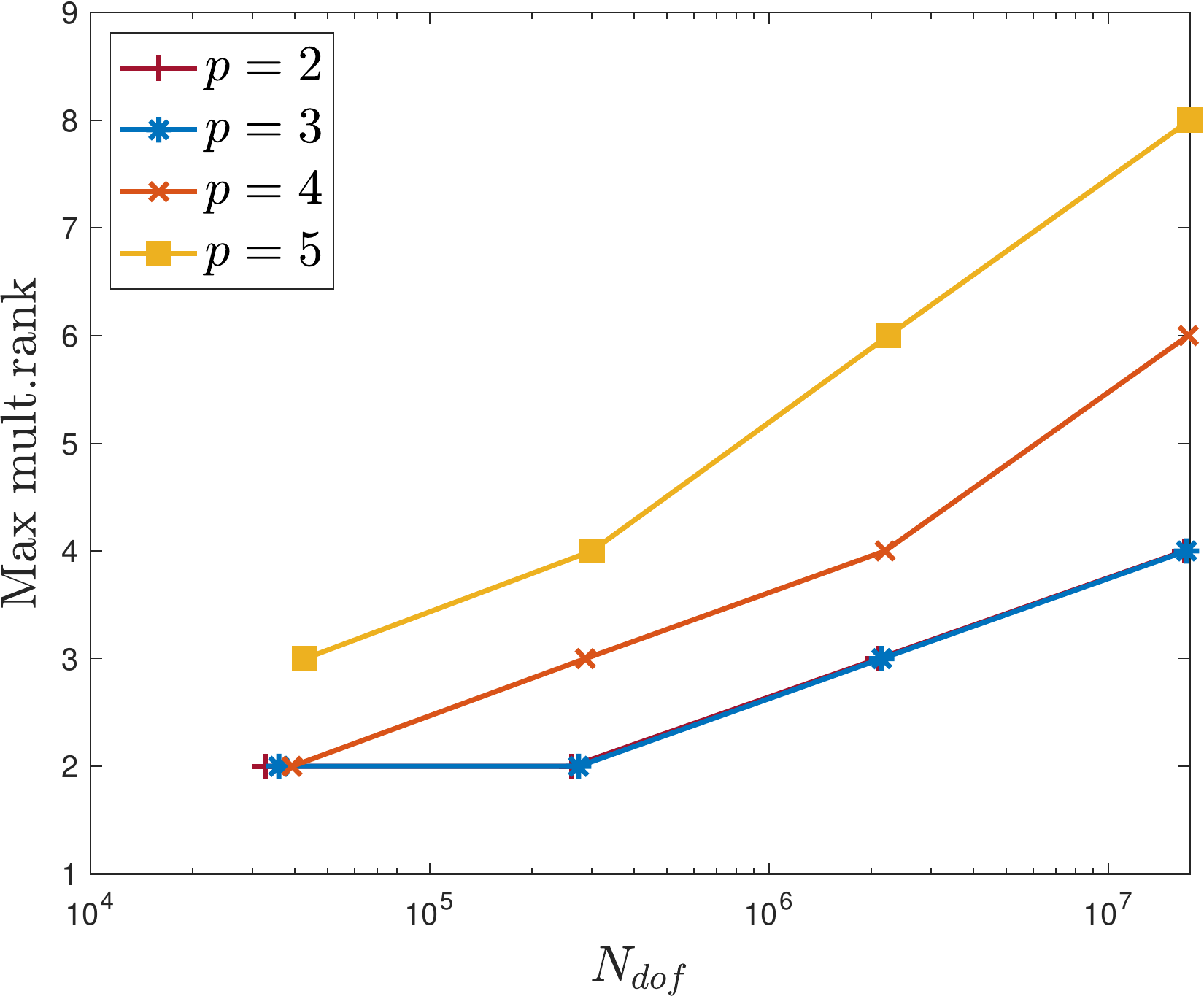}} \quad 
 \subfloat[][Memory compression.\label{fig:mem_TQ_error}]
   {\includegraphics[scale=0.45]{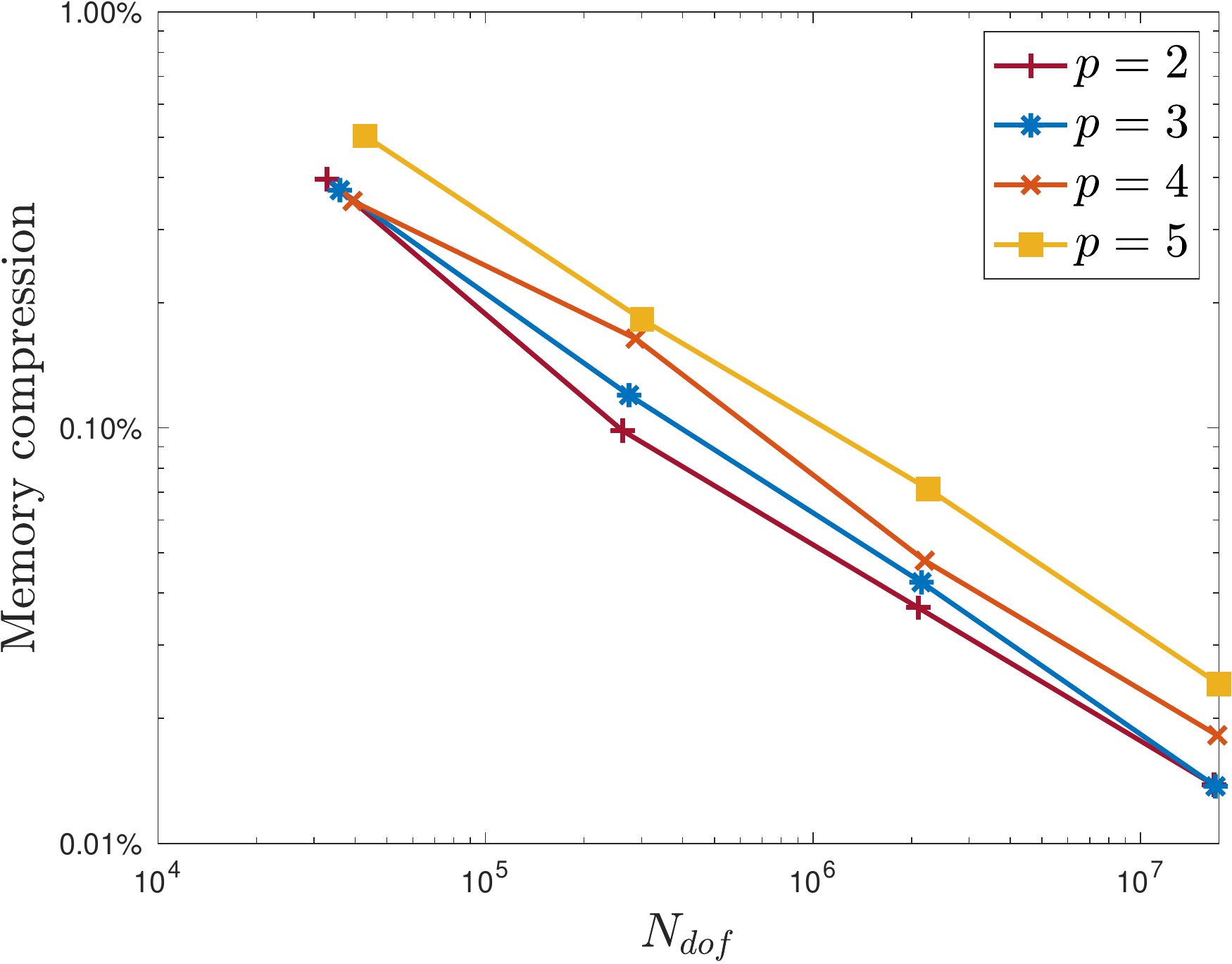}}\quad
 \caption{Results for the thick quarter with variable TPCG tolerance.}
 \label{fig:TQ_results_error}
\end{figure}

For each fixed pair of $h$ and $p$, we estimate  the $L^2$ error of the Galerkin exact solution, and then set the TPCG tolerance $tol$ as this value reduced by a factor of 100. 
This is done in order to 
balance the discretization error and the error introduced by the inexact solution of the linear system, and allows us to observe the right order of convergence for the computed solution. Indeed, if the TPCG tolerance was too loose, then we would observe stagnation in the error convergence curve. On the other hand, if the tolerance was too tight we would perform unuseful iterations. In Figure \ref{fig:quarter_error} we report the  $L^2$ and $H^1$ errors for $p = 2,3,4,5$ and $l = 5,6,7,8$: the orders of convergence exhibit an optimal behavior.
We also report in Figure \ref{fig:rank_TQ_error} the maximum of the multilinear rank of the solution and in Figure \ref{fig:mem_TQ_error} the memory compression. The latter is almost constant with respect to $p$ and decreases as $N_{dof}$ grows, while the maximum of the multilinear rank of the solution mildly grows with respect to $p$ and $N_{dof}$, but it is always relatively small.

In the second test for this geometry, we fix the relative tolerance of TPCG equal  to $10^{-6}$ in order to assess the robustness of our preconditioner. We consider degrees $p=2,3,4,5$ and discretization level $l=7,8,9,10$. 
Note that here we consider finer discretization levels than in the previous test. This is because in the present test we do not compute the $L^2$ and $H^1$ errors, which was the most memory consuming step of the previous test.
The number of iterations, reported in Table \ref{tab:its_TQ_quarter_annulus}, is constant with respect to the number of elements and the degree, assessing the effectiveness of our preconditioning strategy.
 
\begin{figure}[H]
 \centering 
   \subfloat[][Maximum of the multilinear rank of the solution.\label{fig:rank_TQ}] 
   {\includegraphics[scale=0.45]{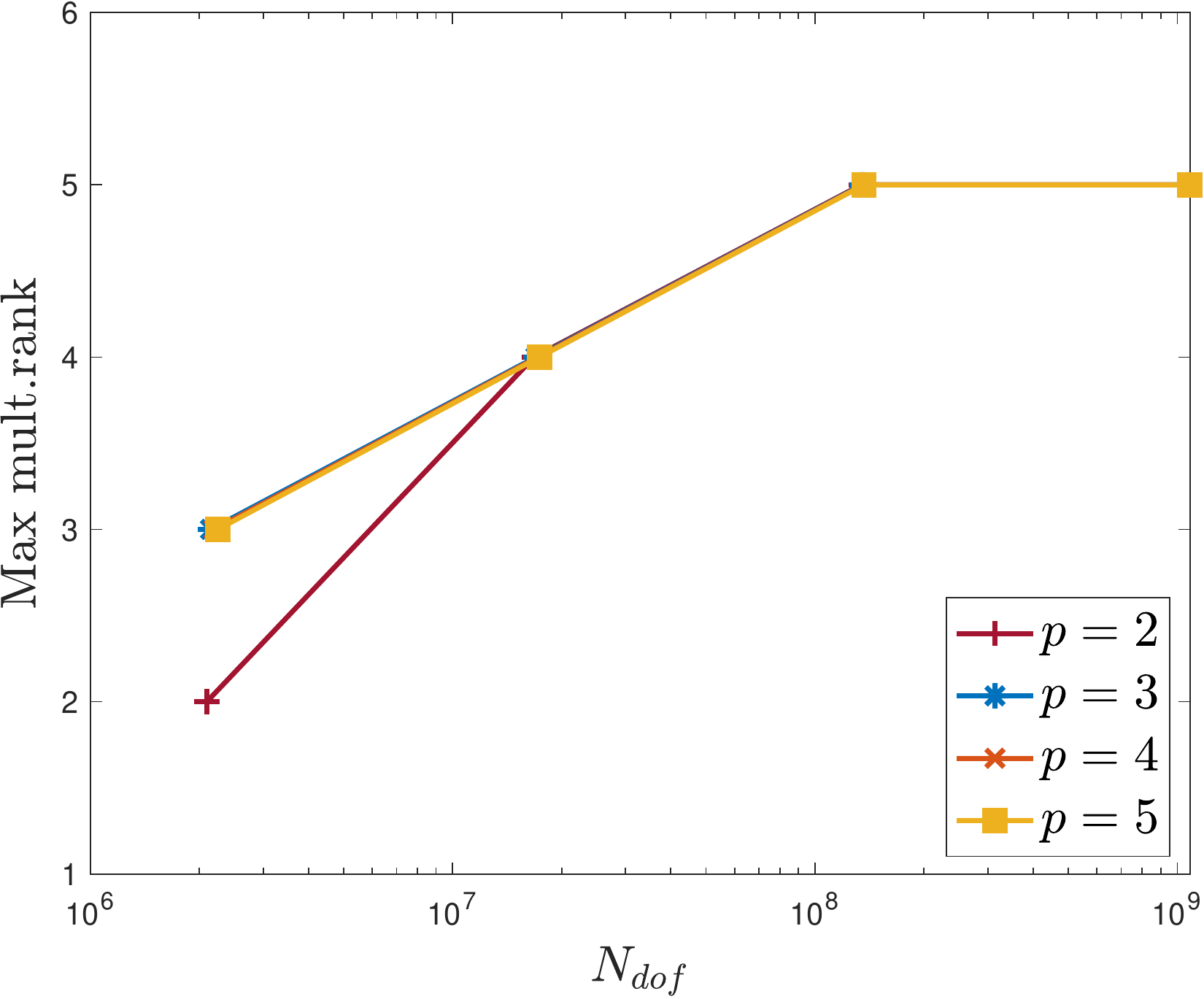}} \quad 
 \subfloat[][Memory compression.\label{fig:mem_TQ}]
   {\includegraphics[scale=0.45]{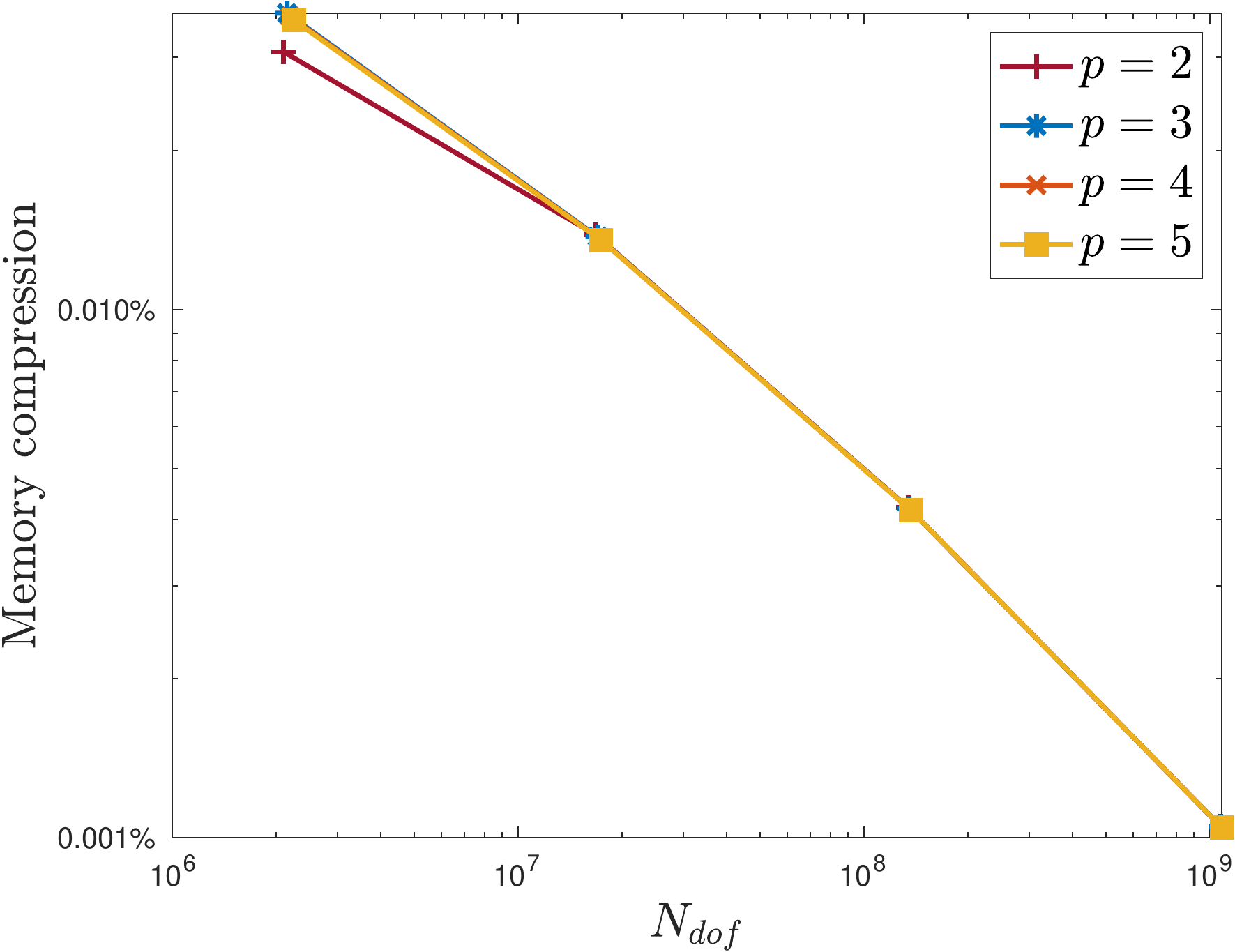}}
 \caption{Results for the thick quarter with $tol=10^{-6}$.}
 \label{fig:TQ_results}
\end{figure}

 {\renewcommand\arraystretch{1.4} 
\begin{table}[H]
\begin{center}
\begin{tabular}{|c|c|c|c|c|}
\hline
 & \multicolumn{4}{|c|}{ \ Iteration number} \\
 \hline
$n_{el}$ & $p=2$  & $p=3$ &$p=4$ & $p=5$  \\
\hline 
128 &   12 &   12  & 12  & 12 \\
\hline
256 &    12 &   12  & 12  & 12 \\ 
\hline
512 &    12 &   12  & 12  & 12\\ 
\hline
1024 &      12 &   12  & 12  & 12  \\ 
\hline
\end{tabular}
\caption{Number of iterations for the thick quarter domain with TPCG
  $tol=10^{-6}$.}
\label{tab:its_TQ_quarter_annulus}
\end{center}
\end{table}}
The maximum of the multilinear rank of the solution is represented in Figure  \ref{fig:rank_TQ}. This value is small, almost independent of $p$, and seems to reach a plateau for the finer discretization levels.

In Figure  \ref{fig:mem_TQ} we report the memory compression. As before, the memory storage is hugely reduced with our low-rank strategy. Moreover, this memory reduction is almost independent on $p$.

\subsection{Spherical shell}
\label{sec:ss_test}

\begin{figure}
\begin{center}
\includegraphics[scale=0.45]{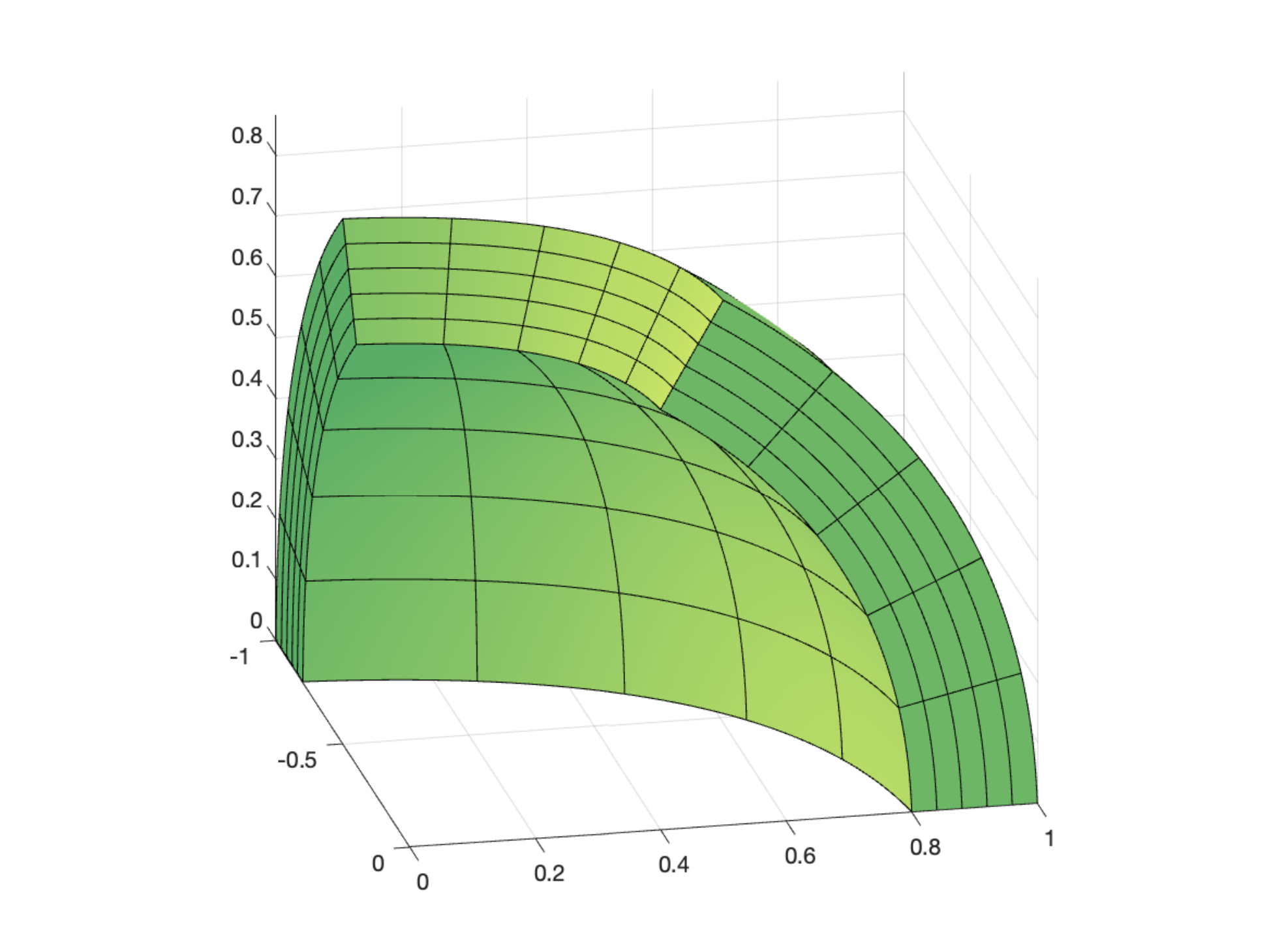}
\caption{Spherical shell domain.}
\label{fig:shell}
\end{center}
\end{figure}

In this test we consider as computational domain a spherical shell, represented in Figure \ref{fig:shell}.

We use $f=1$ as load function, we impose homogeneous Dirichlet boundary condition and we fix the relative tolerance of TPCG as $ 10^{-6}$. We consider degrees $p=2,3,4,5$ and level of discretization $l=7,8,9,10$.

The multilinear ranks of the coefficients of the  stiffness matrix in equation \eqref{eq:matrix_CP}, found using the technique described in \ref{sec:appendix2}, are
$$\begin{bmatrix}
(R_1^{(1,1)},R_2^{(1,1)},R_3^{(1,1)}) &  (R_1^{(1,2)},R_2^{(1,2)},R_3^{(1,2)}) &  (R_1^{(1,3)},R_2^{(1,3)},R_3^{(1,3)}) \\
 (R_1^{(2,1)},R_2^{(2,1)},R_3^{(2,1)}) & (R_1^{(2,2)},R_2^{(2,2)},R_3^{(2,2)}) & ( R_1^{(2,3)},R_2^{(2,3)
},R_3^{(2,3)})\\
(R_1^{(3,1)},R_2^{(3,1)},R_3^{(3,1)}) & (R_1^{(3,2)},R_2^{(3,2)},R_3^{(3,2)}) & (R_1^{(3,3)},R_2^{(3,3)},R_3^{(3,3)})
\end{bmatrix}=
\begin{bmatrix}
(1,1,1) &  (1,1,1) &  (1,1,1) \\
 (1,1,1) & (2,2,1) & (2,2,1)\\
 (1,1,1) &  (2,2,1) & (2,2,1)
\end{bmatrix},
$$
and thus, referring to equation \eqref{eq:tuck_matrix}, we have
$$
  \left(R_1^{A}, R_2^A, R_3^A\right) = (13,13,9).
$$
Note that here the multilinear ranks are higher than for the simpler thick annulus domain considered previously. The number of TPCG iterations, reported in  Table \ref{tab:its_TQ_spherical_shell}, is also higher than in the previous case. This is expected, since here the geometry is not tensor product and as a consequence the preconditioner is less effective. Nevertheless, the number of iterations is almost independent on $p$ and $n_{el}$, showing again that the preconditioner is robust with respect to these parameters.
 

 {\renewcommand\arraystretch{1.4} 
\begin{table}[H]
\begin{center}
\begin{tabular}{|c|c|c|c|c|}
\hline
 & \multicolumn{4}{|c|}{ \ Iteration number} \\
 \hline$n_{el}$ & $p=2$  & $p=3$ &$p=4$  & $p=5$ \\
\hline 
128 &   74   & \z91  & \z95    &  \z88\\
\hline
256 &   79      &  \z94  &  \z99  &  \z96 \\ 
\hline
512 &    85     &  101  &  106  & \z99 \\ 
\hline
1024 &   88      &  101  & 102   & 101  \\ 
\hline
\end{tabular}
\caption{Number of iterations for the spherical shell domain with TPCG $tol=10^{-6}$.}
\label{tab:its_TQ_spherical_shell}
\end{center}
\end{table}}
 
\begin{figure}[H]
 \centering  
     \subfloat[][Maximum of the multilinear rank of the solution.\label{fig:rank_SS}] 
   {\includegraphics[scale=0.45]{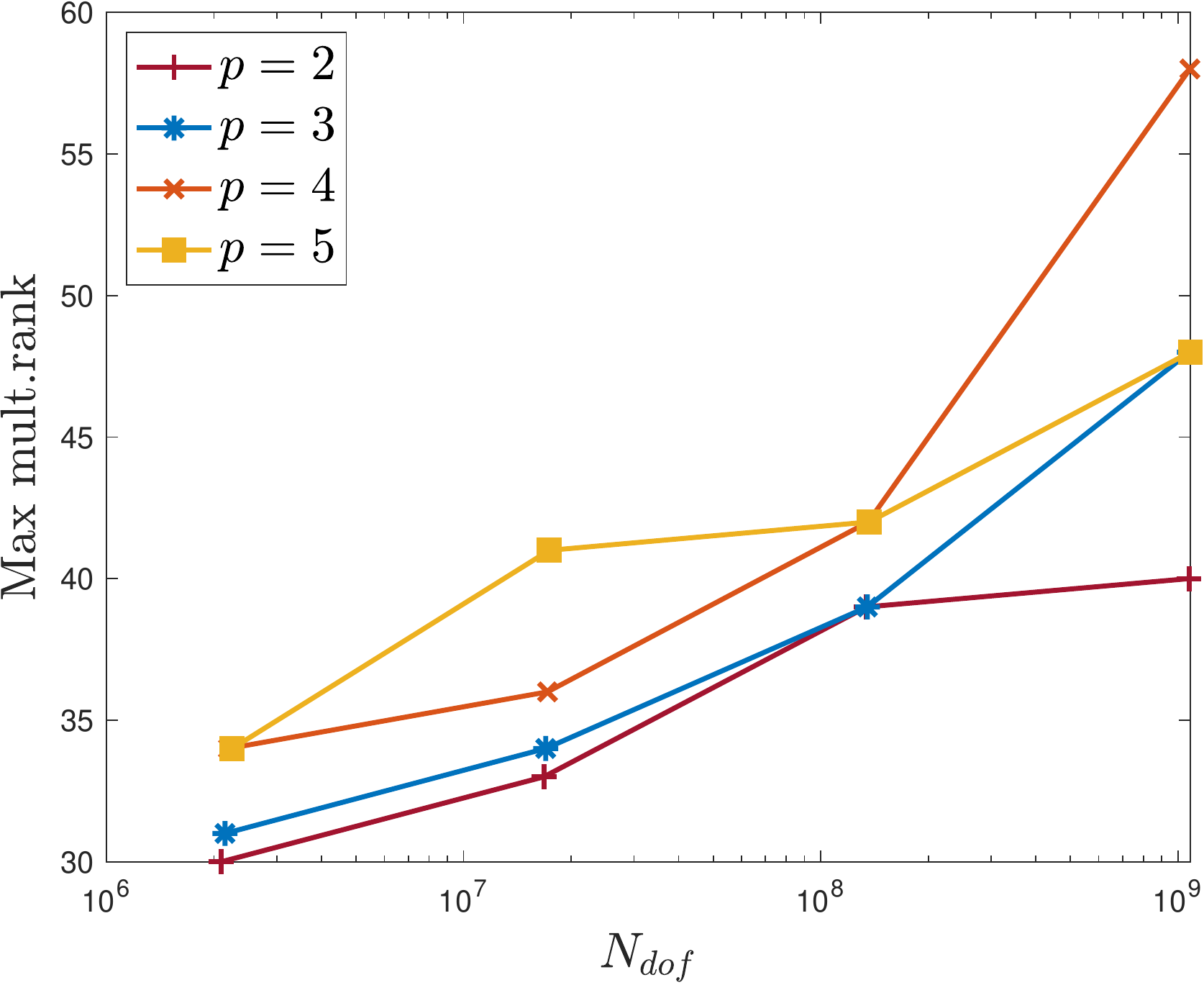}}\quad
 \subfloat[][Memory compression.\label{fig:mem_SS}]
   {\includegraphics[scale=0.45]{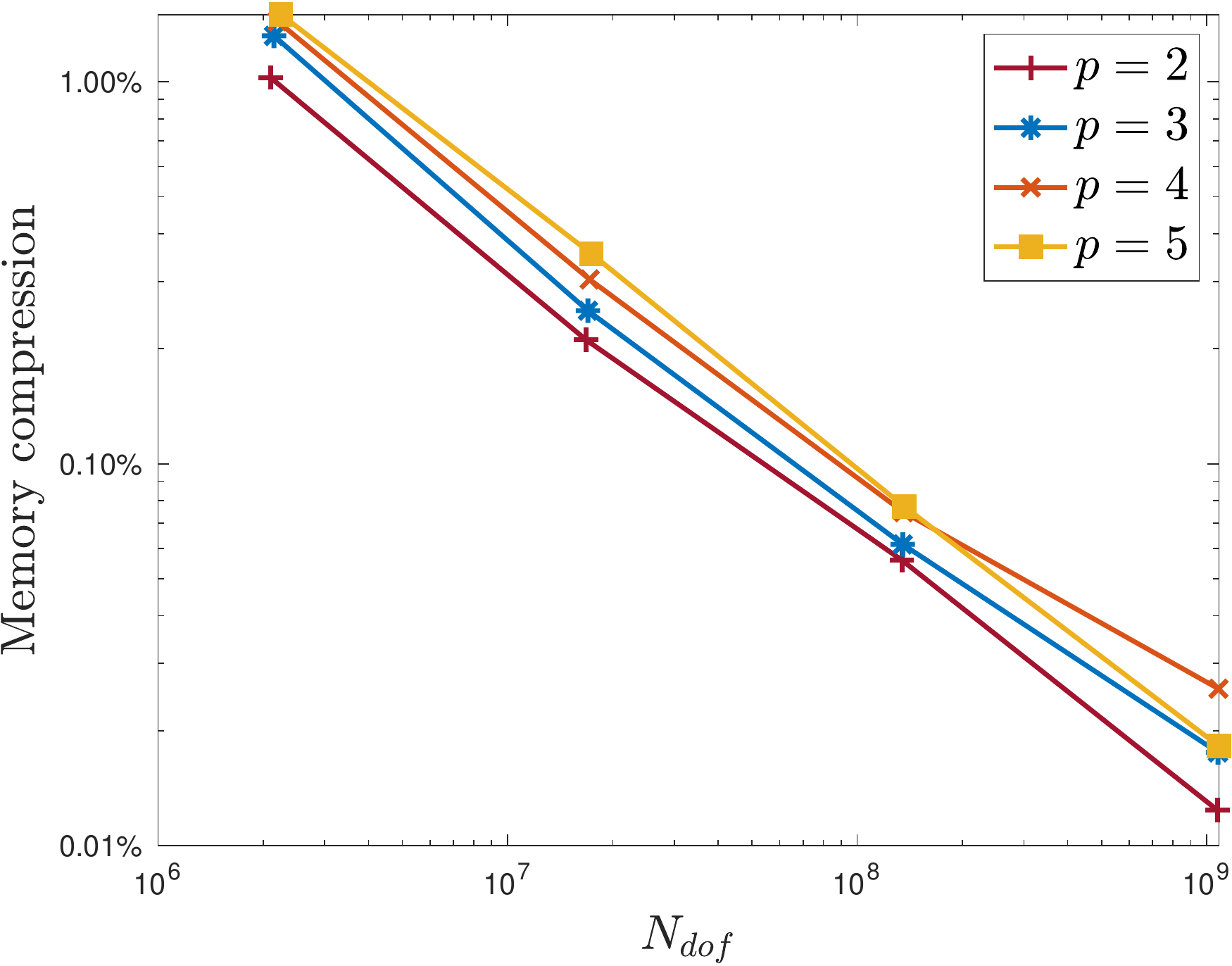}}
   \caption{Results for the spherical shell with TPCG $tol=10^{-6}$.}
 \label{fig:SS_results}
\end{figure}

The maximum of the multilinear rank of the computed solution is represented in Figure \ref{fig:rank_SS}. 
Compared to the previous case, here the rank is 
higher and it also tends to grow as the number of elements and the degree are increased, corresponding to an enrichment of the Galerkin approximation space. This probably indicates that the unknown solution is less suited for a low-rank approximation. 
 
Nevertheless, the memory storage of the low-rank solution with respect to the full numerical solution is hugely reduced (see Figure \ref{fig:mem_SS}). Moreover, similarly to the previous case, this memory reduction is almost independent on $p$. 
 
\begin{rmk}
To assess the quality of the approximated eigendecomposition (see  \ref{sec:appendix} for details) exploited in the preconditioner, we also performed numerical tests using the exact eigenvalues and eigenvector. The numbers of iterations that we obtained do not deviate significantly from the ones shown in Tables \ref{tab:its_TQ_quarter_annulus} and \ref{tab:its_TQ_spherical_shell}.
\end{rmk}

\subsection{Compressible linear elasticity}
\label{sec:le}
  In this test case we consider the compressible linear elasticity problem.
Let $\Omega\subset\mathbb{R}^3$ be the computational domain and let $\partial \Omega$ denote its  boundary. Suppose that  $\partial \Omega=\partial\Omega_D\cup\partial\Omega_N$ with $\partial\Omega_D\cap \partial\Omega_N=\emptyset$ and where $\partial\Omega_D$ has positive measure.   
Let $H^1_D(\Omega):= \left\{v \in H^1(\Omega) \ \middle| \ v=0 \text{ on } \partial\Omega_D \right\}$.  Then, given  $\underline{f} \in  [{L}^2(\Omega)]^3 $ and  $ \ \underline{g} \in   [{L}^2(\partial \Omega_N)]^3$, the  variational formulation of the compressible linear elasticity problem reads: find  $\underline{u}\in [H^1_D(\Omega)]^3$ such that for all   $\underline{v}\in [H^1_D(\Omega)]^3$ 
\begin{equation*}   
b(\underline{u},\underline{v})= \vect{G}( \underline{v}),
\end{equation*}
where we define
\begin{align}
\label{eq:bil_lin}
b(\underline{u},\underline{v}) :=  2 \mu \int_{\Omega} \varepsilon(\underline{u}) : \varepsilon(\underline{v}) \d{\underline{x}} + \lambda \int_{\Omega} \left( \nabla \cdot \underline{u} \right) \left(\nabla \cdot \underline{v}\right) \d{\underline{x}}    , \quad
 \vect{G}( \underline{v}) :=  \int_{\Omega} \underline{f} \cdot \underline{v} \d{\underline{x}} + \int_{\partial \Omega_N} \underline{g}\cdot \underline{v} \ds, 
\end{align}
$ {\varepsilon}(\underline{v}) := \frac{1}{2} \left(\nabla \underline{v}+ (\nabla \underline{v})^T\right)$ is the symmetric gradient,   while   $\lambda$ and $\mu$ denote the material Lam\'{e} coefficients.  
 We consider the isogeometric discretization presented in \cite[Section~6.2]{Beirao2014}. With a proper ordering of the degrees of freedom, the resulting linear system has a $3 \times 3$ block structure as a consequence of the vectorial nature of the isogeometric space (see \cite{Beirao2014}). Thus, we have to solve $ \mathbf{A}\vect{\underline{x}}=\vect{\underline{f}} $
where 
\begin{equation}
  \label{eq:block-system-elasticity}
  \mathbf{A}:=\begin{bmatrix}
 \mathbf{A}_{1,1} &  \mathbf{A}_{1,2} &  \mathbf{A}_{1,3}\\
  \mathbf{A}_{2,1} &  \mathbf{A}_{2,2} &  \mathbf{A}_{2,3}\\
   \mathbf{A}_{3,1} &  \mathbf{A}_{3,2} &  \mathbf{A}_{3,3}
 \end{bmatrix},\quad \vect{\underline{x}}:=\begin{bmatrix}
 \vect{x}_1\\ \vect{x}_2\\ \vect{x}_3
 \end{bmatrix},\quad \text{and} \quad \vect{\underline{f}}:=\begin{bmatrix}
 \vect{f}_1\\ \vect{f}_2\\ \vect{f}_3
 \end{bmatrix}.
\end{equation}
We look for a solution $\widetilde{\vect{\underline{x}}}=\begin{bmatrix}
 \widetilde{\vect{x}}_1\\ \widetilde{\vect{x}}_2\\ \widetilde{\vect{x}}_3
 \end{bmatrix}$ such that each of its component $\widetilde{\vect{x}}_1, \widetilde{\vect{x}}_2$ and $\widetilde{\vect{x}}_3$ is a  Tucker vector. In order to have a system matrix and a right-hand side with a compatible structure, we apply the Tucker approximation strategy described in \ref{sec:appendix2} independently to each block of the system matrix and of the right-hand side. In other words, each $\mathbf{A}_{i,j}$ is approximated by a Tucker matrix and each $\vect{f}_j$ is approximated by a Tucker vector for $i,j=1,2,3$, up to a fixed relative tolerance.

 Following \cite{Bosy2020}, as a preconditioner one can  take a block
diagonal approximation of $\mathbf{A}$ (seen as a $3\times 3$ block
matrix, as in \eqref{eq:block-system-elasticity}) where each of the three blocks
is computed in parametric coordinates, that is, without dependence on the
geometry parameterization (see \cite[Lemma 1]{Bosy2020}). Here, each
diagonal block is further approximated  following the same steps
described in Section \ref{sec:preconditioner}.  For all diagonal blocks we choose a relative tolerance $\epsilon_{rel}^{prec}=10^{-1}$, as in the Poisson tests. The ranks $R_P$ that we found are comparable, though slightly higher, to those reported in Table \ref{tab:ranks_table} for the Poisson problem. 

Matrix-vector products, scalar products, truncations and sums are handled block-wise using the same techniques described in Section \ref{sec:tensor_calculus} and Section \ref{sec:truncation}. 

\begin{figure}
\centering
  \subfloat[][Column domain.\label{fig:def_col}] 
   {\includegraphics[scale=0.45]{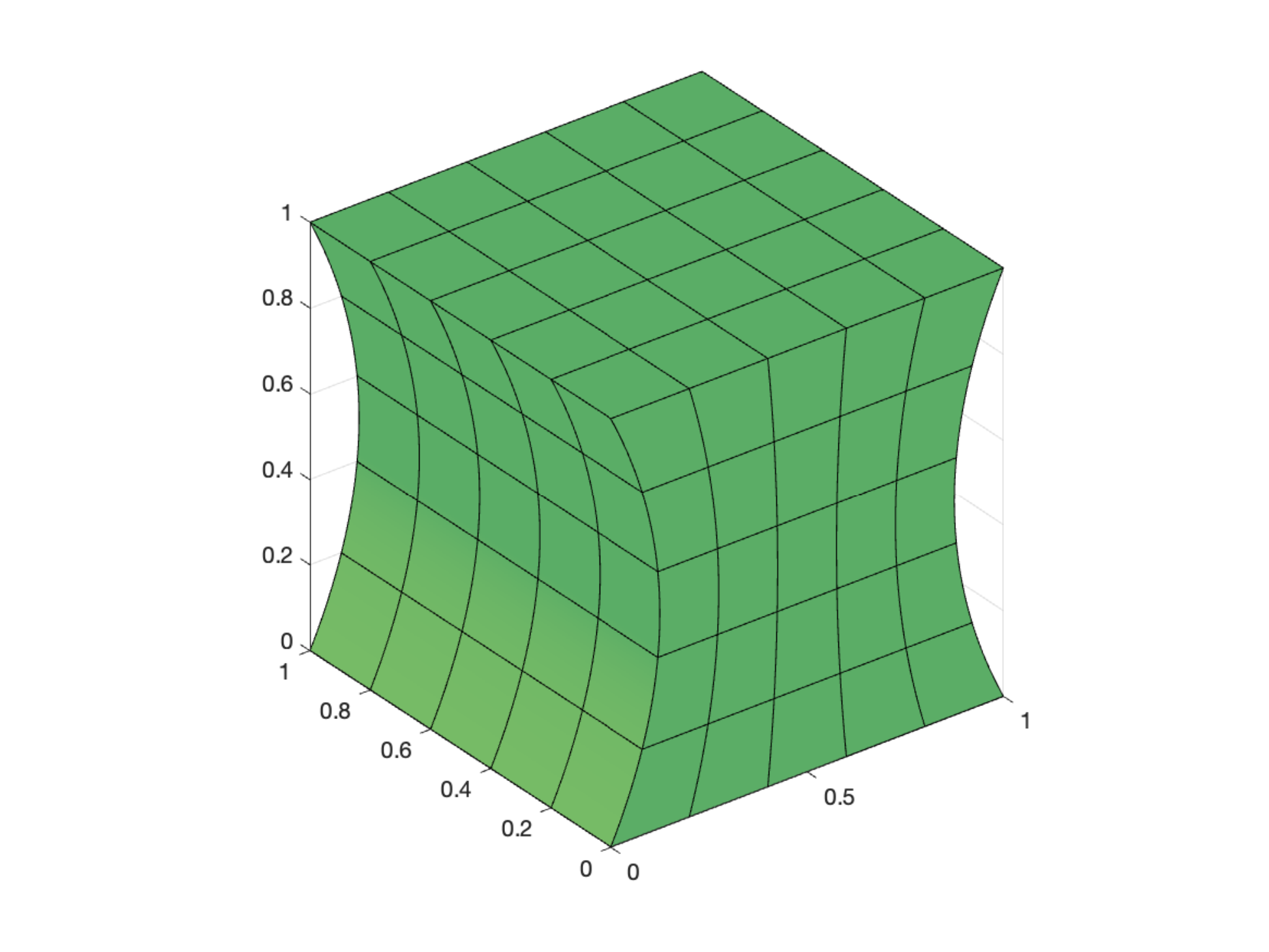}} 
 \subfloat[][Numerical solution.\label{fig:num_sol_LE}]
  {\includegraphics[scale=0.40]{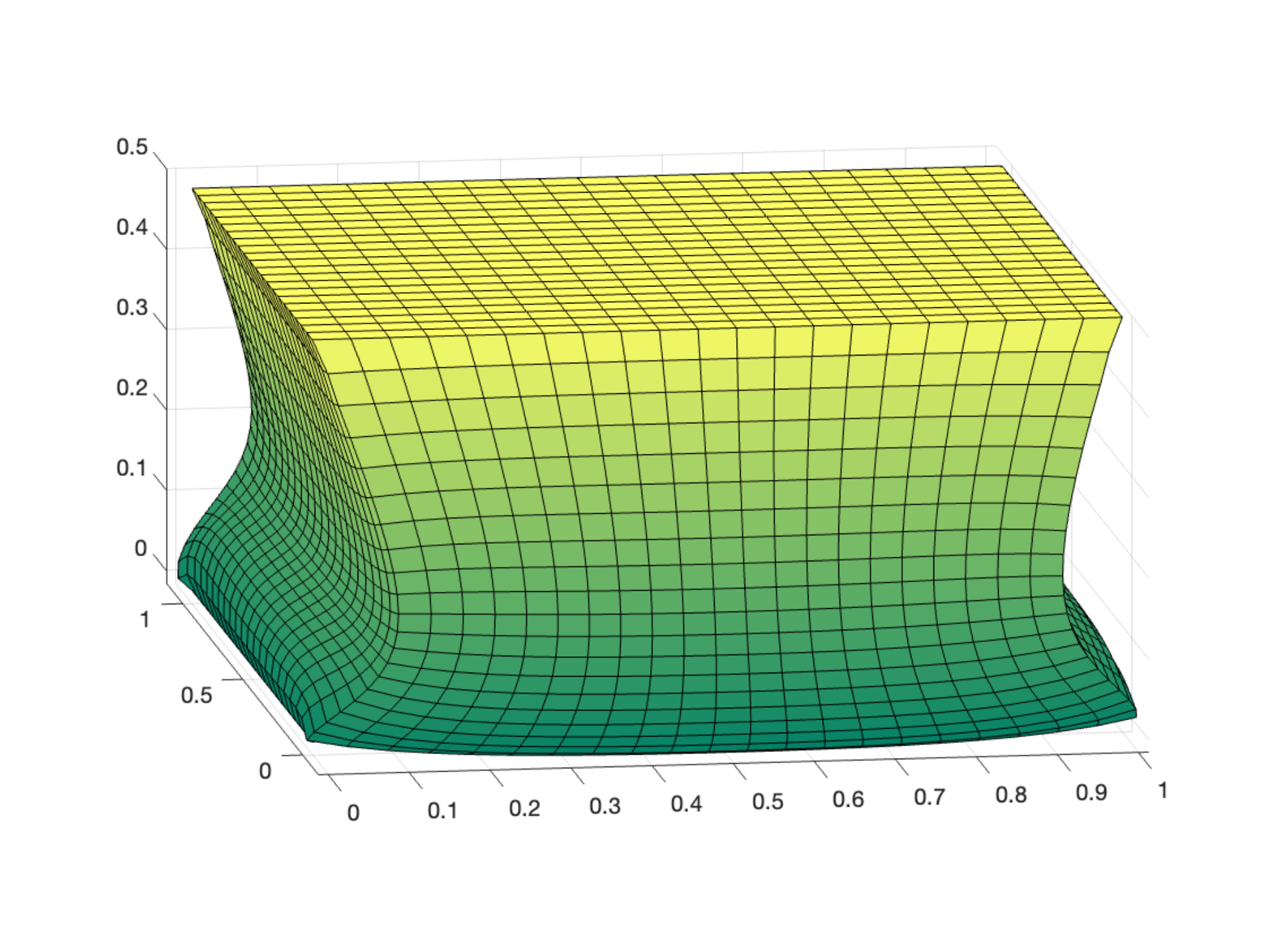}}
  \caption{Linear elasticity: initial domain and numerical solution.}
\end{figure}

For the numerical test, we consider as computational domain a column with square section and with two faces described by a quadratic function, represented in Figure \ref{fig:def_col}. 

We choose $\underline{f}=[0, 0, -1]^T$, homogeneous Dirichlet boundary condition on the face $\{z=0\}$,  we impose $\underline{u}= [0, 0, -0.5]^T$ on the face $\{z=1\}$   and homogeneous Neumann boundary condition on the other sides. The values of the  Lam\'{e} parameters are set equal to  $\lambda=\frac{0.3}{0.52}$ and $\mu =\frac{1}{2.6}$, corresponding to the choice $E=1$ for the Young Modulus and $\nu=0.3$ for the Poisson's ratio  (recall that $\lambda := \frac{E\nu}{(1+\nu)(1-2\nu)}$ and $\mu:=\frac{E}{2(1+\nu)}$).  We fix the TPCG tolerance equal to $10^{-6}$. We report the numerical solution that we obtained in Figure \ref{fig:num_sol_LE}.

Referring to equation \eqref{eq:tuck_matrix},  we have that the multilinear ranks of the approximated matrix blocks are
$$
  \left(R_1^{A_{1,1}}, R_2^{A_{1,1}}, R_3^{A_{1,1}}\right) =  \left(R_1^{A_{2,2}}, R_2^{A_{2,2}}, R_3^{A_{2,2}}\right) =  \left(R_1^{A_{3,3}}, R_2^{A_{3,3}}, R_3^{A_{3,3}}\right) = (6,5,6).
$$
and 
\begin{align*}
  \left(R_1^{A_{1,2}}, R_2^{A_{1,2}}, R_3^{A_{1,2}}\right) & =  \left(R_1^{A_{2,1}}, R_2^{A_{2,1}}, R_3^{A_{2,1}}\right) =  \left(R_1^{A_{1,3}}, R_2^{A_{1,3}}, R_3^{A_{1,3}}\right) = \left(R_1^{A_{3,1}}, R_2^{A_{3,1}}, R_3^{A_{3,1}}\right) \\
  &=\left(R_1^{A_{2,3}}, R_2^{A_{2,3}}, R_3^{A_{2,3}}\right) = \left(R_1^{A_{3,2}}, R_2^{A_{3,2}}, R_3^{A_{3,2}}\right) = (2,2,2).
\end{align*}

In particular, the multilinear ranks of the coefficients of each block are
$$
\begin{bmatrix}
\begin{bmatrix}
(2,1,2) &  \mathbf{0} &  (1,1,1) \\
 \mathbf{0} & (1,1,1) &  \mathbf{0} \\
 (1,1,1) &  \mathbf{0}  & (1,1,1)
\end{bmatrix} &  \begin{bmatrix}
 \mathbf{0} &   (1,1,1) &  \mathbf{0} \\
  (1,1,1) &  \mathbf{0}&  \mathbf{0} \\
  \mathbf{0}&  \mathbf{0}  &  \mathbf{0}
\end{bmatrix}& \begin{bmatrix}
 \mathbf{0} &  \mathbf{0} &  (1,1,1) \\
 \mathbf{0} &  \mathbf{0}&  \mathbf{0} \\
 (1,1,1) &  \mathbf{0}  &  \mathbf{0}
\end{bmatrix}\\
& & &\\
\begin{bmatrix}
 \mathbf{0} &   (1,1,1) &  \mathbf{0} \\
  (1,1,1) &  \mathbf{0}&  \mathbf{0} \\
  \mathbf{0}&  \mathbf{0}  &  \mathbf{0}
\end{bmatrix}& \begin{bmatrix}
(2,1,2) &  \mathbf{0} &  (1,1,1) \\
 \mathbf{0} & (1,1,1) &  \mathbf{0} \\
 (1,1,1) &  \mathbf{0}  & (1,1,1)
\end{bmatrix} & \begin{bmatrix}
 \mathbf{0} &  \mathbf{0} &  \mathbf{0}\\
 \mathbf{0} &  \mathbf{0}&  (1,1,1)\\
 \mathbf{0}&   (1,1,1) &  \mathbf{0}
\end{bmatrix}\\
& & &\\
\begin{bmatrix}
 \mathbf{0} &  \mathbf{0} &  (1,1,1) \\
 \mathbf{0} &  \mathbf{0}&  \mathbf{0} \\
 (1,1,1) &  \mathbf{0}  &  \mathbf{0}
\end{bmatrix}&\begin{bmatrix}
 \mathbf{0} &  \mathbf{0} &  \mathbf{0}\\
 \mathbf{0} &  \mathbf{0}&  (1,1,1)\\
 \mathbf{0}&   (1,1,1) &  \mathbf{0}
\end{bmatrix} & \begin{bmatrix}
(2,1,2) &  \mathbf{0} &  (1,1,1) \\
 \mathbf{0} & (1,1,1) &  \mathbf{0} \\
 (1,1,1) &  \mathbf{0}  & (1,1,1)
\end{bmatrix}
\end{bmatrix}.
$$

We consider degrees $p=3,4,5$ and discretization levels $l= 7,8,9,10$. The number of iterations, reported in Table \ref{tab:its_LE} is almost independent of $p$ and grows only mildly with the number of elements $n_{el}$.  The only exception is represented by the case $p=3$ and $n_{el}=1024$, where we observe a stagnation in the residual. This issue can be fixed by decreasing the minimum relative tolerance $\epsilon_{\min}$ used for the dynamical truncation. 

 {\renewcommand\arraystretch{1.4} 
\begin{table}[H]
\begin{center}

\begin{tabular}{|c|c|c|c|}
\hline
 & \multicolumn{3}{|c|}{ \ Iteration number} \\
 \hline
$n_{el}$ & $p=3$ &$p=4$  & $p=5$ \\
\hline 
128 &        20  &  19    & 18  \\
\hline
256 &      25         &   20  &  20  \\ 
\hline
512 &     29        &  23     &  22  \\  
\hline
1024 &    43    & 28       & 28    \\
\hline
\end{tabular}

\caption{Number of iterations for the deformed column domain  with TPCG $tol=10^{-6}$.}
\label{tab:its_LE}
\end{center}
\end{table}}

\begin{figure}[H]
 \centering  
     \subfloat[][Maximum of the multilinear ranks of all the components of the solution. \label{fig:rank_LE}] 
   {\includegraphics[scale=0.45]{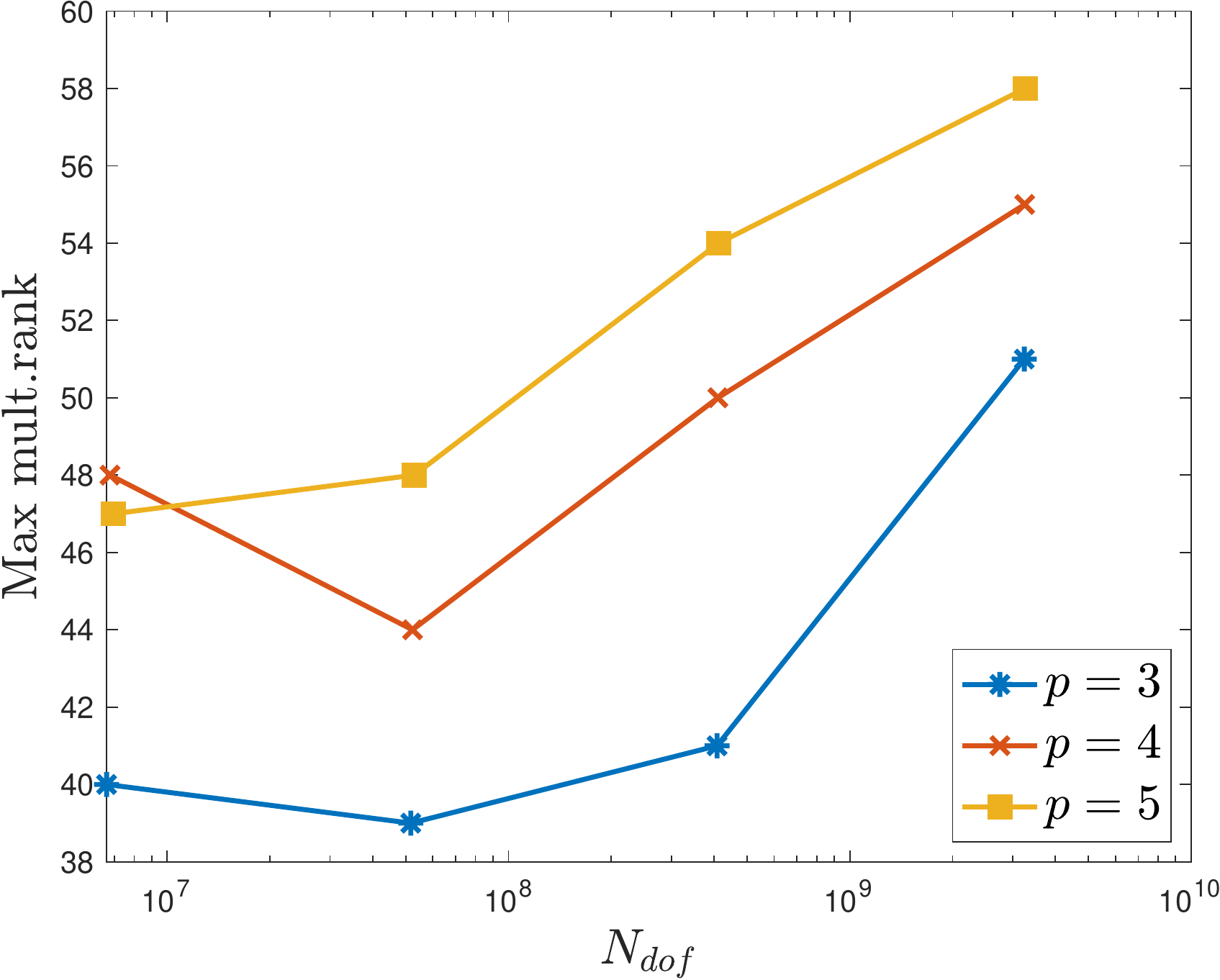}}
   \subfloat[][Memory compression.\label{fig:mem_LE}]
   {\includegraphics[scale=0.45]{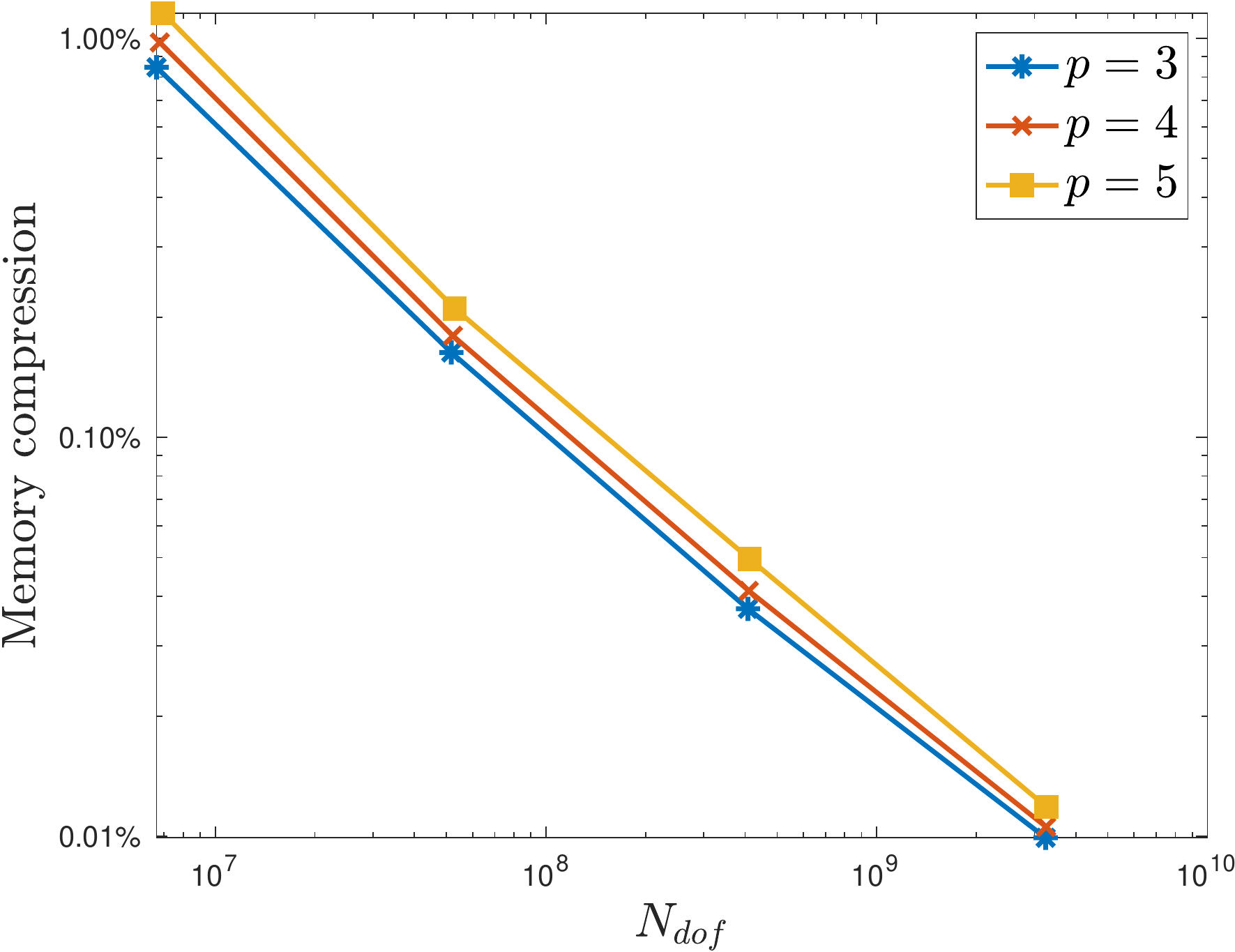}}
   \caption{Results for the liner elasticity problem.}
 \label{fig:LE_results}
\end{figure}

We represent in Figure \ref{fig:rank_LE} the maximum among the maxima multilinear ranks of the three components of the computed solution. 
We observe that, as in the spherical shell test of Section \ref{sec:ss_test}, the ranks grow as  the  degree and number of elements is increased. 
  
If $\widetilde{\vect{x}} =\begin{bmatrix}
 \widetilde{\vect{x}}_1\\ \widetilde{\vect{x}}_2\\ \widetilde{\vect{x}}_3
 \end{bmatrix}$ is the  solution of the linear elasticity problem with $\widetilde{\vect{x}}_i$ Tucker vectors    represented by the Tucker tensors
$\widetilde{\mathcal{X}}^{(i)}=\widetilde{\core{X}}^{(i)}\times_3\widetilde{\factor{X}}^{(i)}_3\times_2\widetilde{\factor{X}}^{(i)}_2\times_1\widetilde{\factor{X}}^{(i)}_1$
with $\widetilde{\core{X}}^{(i)}\in\mathbb{R}^{r^{(i)}_1\times r^{(i)}_2\times r^{(i)}_3}$ and
$\widetilde{\factor{X}}^{(i)}_i\in\mathbb{R}^{n^{(i)}_j\times
  r^{(i)}_j}$ for $i, j=1,2,3$,  we define the memory compression
percentage  for compressible linear elasticity problem as
\begin{equation}
\text{memory compression}:= \frac{\sum_{i=1}^3r^{(i)}_1r^{(i)}_2r^{(i)}_3 + r_1^{(i)}n_1^{(i)}+ r_2^{(i)}n_2^{(i)}+ r_3^{(i)}n_3^{(i)}}{\sum_{i=1}^3 n^{(i)}_1n^{(i)}_2n^{(i)}_3}\cdot 100.
\end{equation}

As in case of Poisson, the memory storage for the solution is hugely reduced with respect to the full solution, as reported in Figure \ref{fig:mem_LE}. Note in particular that the  memory compression is  around  $0.01\%$  for the finest level of discretization.

\section{Conclusions}
\label{sec:conclusions}
In this work we presented a low-rank IgA solver for the Poisson
problems. The stiffness matrix and right-hand side vector are
approximated as   Tucker matrix and  Tucker vector,
respectively. The solution of the linear system is then computed  in
Tucker format by the  Truncated Preconditioned Conjugate Gradient
method. We have designed a novel preconditioner combining the  Fast
Diagonalization method, an approximation of the eigenvectors that allows
multiplication by the  Fast Fourier Transform, and a low-rank approximation of the inverse eigenvalues using a sum of exponential function.
%
The cost of the setup and application of the preconditioner is almost
linear with respect to the number of univariate dofs.

The numerical tests, performed  on Poisson and on linear elasticity
benchmarks, confirm the expected  low memory
storage, almost independent  of the  degree $p$. 

Future research directions will go towards the mathematical analysis of  the eigenvalues/eigenvectors approximation used for the preconditioner and reported in \ref{sec:appendix}.
We also plan to extend the present low-rank strategy to multi-patch
domains.

Moreover, we would like to improve the preconditioning strategy by allowing it to take into account any scalar coefficient appearing in the equation. If the coefficient is continuous, this could  be done by introducing a diagonal scaling similar to the one used in \cite{Loli2022}, combined with a low-rank compression. On the other hand, the case of a discontinuous coefficient is more challenging and requires further investigation.


Finally, we will consider   tensor formats more  efficient 
in higher dimensions, such as  Tensor Train.
 
\section*{Acknowledgments}

The   authors are members of the Gruppo Nazionale Calcolo
Scientifico-Istituto Nazionale di Alta Matematica (GNCS-INDAM), and in
particular the first author was partially supported by INDAM-GNCS Project 2023 ``Efficienza ed analisi di metodi numerici innovativi per la risoluzione di PDE'' and the third  author was partially supported by INDAM-GNCS Project 2022 ``Metodi numerici efficienti e innovativi per la
risoluzione di PDE''.  The second author acknowledges the
contribution of the National Recovery and Resilience Plan, Mission 4
Component 2 - Investment 1.4 - NATIONAL CENTER FOR HPC, BIG DATA AND
QUANTUM COMPUTING, spoke 6. 

\appendix

\section{Low-rank approximation of the linear system}
\label{sec:appendix2} 

We present here the procedure that we use to  find a Tucker approximation of $\vect{f}$ and $\matb{A}$ of the form \eqref{eq:Tucker} and \eqref{eq:tucker_matrix}, respectively.    
%
 First, we write the system matrix and the right-hand side as
 \begin{equation}
 \label{eq:matrix_param}
 [\matb{A}]_{i,j}=\int_{\widehat{\Omega}}\left(\nabla\widehat{B}_i(\vect{\eta})\right)^T\matb{Q(\vect{\eta})}\nabla\widehat{B}_j(\vect{\eta})\d\widehat{\Omega}=\sum_{k,l=1}^d \int_{\widehat{\Omega}} [\matb{Q}(\vect{\eta})]_{k,l}\frac{\partial}{\partial\eta_l}\widehat{B}_i(\vect{\eta})  \frac{\partial}{\partial\eta_k}\widehat{B}_j(\vect{\eta})\d\widehat{\Omega}
\end{equation}
for $i,j=1,\dots,N_{dof}$ and
 \begin{equation}
 \label{eq:rhs_param}
 [\vect{f}]_{i}=\int_{\widehat{\Omega}} \widehat{B}_i(\vect{\eta})\omega(\vect{\eta})\d\widehat{\Omega},
\end{equation}
 for $i=1,\dots,N_{dof}$, respectively, where
 $$
\matb{Q}(\vect{\eta}):= \text{det}(J_{\vect{F}}(\vect{\eta}))J_{\vect{F}}^{-1}(\vect{\eta})J_{\vect{F}}(\vect{\eta})^{-T}\in\mathbb{R}^{d\times d},\quad \omega(\vect{\eta}):= \text{det}(J_{\vect{F}}(\vect{\eta})) f(\vect{\eta})\in\mathbb{R}
 $$
 and $J_{\matb{F}}$ is the Jacobian of $\vect{F}$.
If the entries of $\matb{Q}$ and $\omega$ are sums of separable functions, then \eqref{eq:matrix_param}-\eqref{eq:rhs_param} can be written as the sum of the product of univariate integrals from which expressions of the form \eqref{eq:Tucker}-\eqref{eq:tucker_matrix} immediately follow. However, for many non-trivial geometries and load functions, this is not the case.

 We focus on $d=3$, that is our case of interest, and   we consider the strategy proposed in \cite{Dolgov2021} and implemented in the Chebfun Matlab toolbox \cite{Driscoll2014}.  For dimensions higher than 3, one could use the approach described in \cite{Mantzaflaris2017}.
Fixed a relative tolerance $\epsilon$, we apply \textsf{chebfun3f} function \cite{Dolgov2021} to the entries of $\matb{Q}(\boldsymbol{\eta})$ and to $\omega(\boldsymbol{\eta})$ and   we find separable approximations
 $$
[\matb{Q}(\boldsymbol{\eta})]_{k,l}\approx [\widetilde{\matb{Q}}(\boldsymbol{\eta})]_{k,l}:=  \sum_{i_1=1}^{p^{(k,l)}_1}\sum_{i_2=1}^{p^{(k,l)}_2}\sum_{i_3=1}^{p^{(k,l)}_3}[\mathcal{Q}^{(k,l)}]_{i_1,i_2,i_3}T_{i_1}(\eta_1)T_{i_2}(\eta_2) T_{i_3}(\eta_3) \quad \text{ for } k,l=1,2,3 $$
and 
$$ \omega(\boldsymbol{\eta})\approx \widetilde{\omega}(\boldsymbol{\eta}) :=  \sum_{i_1=1}^{q_1}\sum_{i_3=1}^{q_2}\sum_{i_3=1}^{q_3}[\mathcal{W}]_{i_1,i_2,i_3}T_{i_1}(\eta_1) T_{i_2}(\eta_2)  T_{i_3}(\eta_3),
$$
 where $T_k(\eta_k)=\cos((k-1)\cos^{-1}(\eta_k))$ is the $k$-th Chebyschev polynomial, and $\mathcal{Q}^{(k,l)}=\core{Q}^{(k,l)} \times_3 \matb{Q}_{3}^{(k,l)}  \times_2  \matb{Q}_{2}^{(k,l)} \times_1  \matb{Q}_{1}^{(k,l)}$ and $\mathcal{W} :=\core{W} \times_3 \matb{W}_{3}  \times_2  \matb{W}_{2} \times_1  \matb{W}_{1}$ are the  tensors in Tucker format. Here $\core{Q}^{(k,l)}\in\mathbb{R}^{R^{(k,l)}_1\times R^{(k,l)}_2\times R^{(k,l)}_3}$ and  $\core{W}\in\mathbb{R}^{R_1\times R_3\times R_3}$  are the core tensors and $\matb{Q}^{(k,l)}_{t}\in\mathbb{R}^{p_t^{(k,l)}\times R^{(k,l)}_t}$ and $\matb{W}_{t}\in\mathbb{R}^{q_t\times R_t} $ are factor matrices of $ {\mathcal{Q}}^{(k,l)}$ and ${\mathcal{W}}$, respectively.  These approximations satisfy
$$
|[\matb{Q}(\boldsymbol{\eta})]_{k,l}-[\widetilde{\matb{Q}}(\boldsymbol{\eta})]_{k,l}|\leq 10 \epsilon \quad \text{and}\quad |\omega(\boldsymbol{\eta})-\widetilde{\omega}(\boldsymbol{\eta})|\leq 10 \epsilon
$$
at the Halton points \cite{Niederreiter1992}.  The   algorithm behind \textsf{chebfun3f}  combines a tensorized Chebyshev interpolation and a low-rank Tucker approximation of the evaluation tensor, which is never computed in full format. We refer to the original paper \cite{Dolgov2021}  for the details.
We remark that, as  $\mathbf{Q}(\eta)$  is symmetric, we can further reduce the storage and computational costs, considering only the approximation of $6$ entries of $\mathbf{Q}(\eta)$.

We can write 
\begin{equation}
\label{eq:approxQ}
 [\widetilde{\matb{Q}}(\vect{\eta})]_{k,l}=  \left( {\vect{T}}_{{p_3}^{(k,l)}}(\eta_3)\matb{Q}_3^{(k,l)}\otimes  {\vect{T}}_{p_2^{(k,l)}}(\eta_2)\matb{Q}_2^{(k,l)} \otimes  {\vect{T}}_{p_1^{(k,l)}}(\eta_1)\matb{Q}_1^{(k,l)} \right) \vec(\core{Q}^{(k,l)} )
\end{equation}
and
\begin{equation}
\label{eq:approxW}  \widetilde{\omega}(\vect{\eta})=  \left({\vect{T}}_{q_3}(\eta_3)\matb{W}_3\otimes {\vect{T}}_{q_2}(\eta_2)\matb{W}_2 \otimes {\vect{T}}_{q_1}(\eta_1)\matb{W}_1 \right) \vec(\core{W}),
\end{equation}
where    ${\vect{T}}_m(x):=[T_1(x), \dots , T_m(x) ]\in\mathbb{R}^{m} $ is a vector that collects the Chebyschev polynomials up to the $m$-th one.
Replacing $[\matb{Q}(\vect{\eta})]_{k,l}$ by its approximation \eqref{eq:approxQ} in \eqref{eq:matrix_param} and associating  the indexes $i$ and $j$ with the corresponding indexes of the univariate functions, i.e. $i\leftrightarrow (i_1,i_2,i_3)$ and $j\leftrightarrow (j_1,j_2,j_3)$,  we get
 \begin{align*} 
[\matb{A}]_{i,j} \approx[\widetilde{\matb{A}}]_{i,j}&:=\sum_{k,l=1}^3   \int_{\widehat{\Omega}}   [\widetilde{\matb{Q}}(\vect{\eta})]_{k,l} 
\frac{\partial}{\partial\eta_l}\widehat{B}_i(\vect{\eta})  \frac{\partial}{\partial\eta_k}\widehat{B}_j(\vect{\eta})\d\widehat{\Omega}\\ 
&=\sum_{k,l=1}^3 \sum_{r_3=1}^{R^{(k,l)}_3} \sum_{r_2=1}^{R^{(k,l)}_2} \sum_{r_1=1}^{R^{(k,l)}_1}[\core{Q}^{(k,l)}]_{r_1,r_2,r_3}   \prod_{t=1}^3 \int_0^1 \delta^{(l,t)}(\widehat{b}_{i_t}(\eta_t) )\delta^{(k,t)}(\widehat{b}_{j_t}(\eta_t))\matb{Q}^{(k,l)}_{t}(:,r_t) \cdot{\vect{T}}_{p_t^{(k,l)}}(\eta_t) \d\eta_t,
 \end{align*}
where $\matb{Q}^{(k,l)}_{t}(:,r_t)\in\mathbb{R}^{m_t}$ denotes the $r_t$-th column of $\matb{V}_t^{(k,l)}$ and $\delta^{(l,t)}$ is defined as the operator that acts on a function $f$ as
$$
\delta^{(l,t)}(f):=\left\{\begin{array}{l}
\frac{\partial f}{\partial \eta_t} \quad \text{if}\quad t=l,\\
f\quad \quad \text{otherwise}.
\end{array} \right.
$$
Finally, we have that $\widetilde{\matb{A}}$ is the sum of Tucker matrices of the form \eqref{eq:tucker_matrix}
\begin{equation}
\label{eq:matrix_CP}
\widetilde{\matb{A}} = \sum_{k,l=1}^3\sum_{r_3=1}^{R^{(k,l)}_3} \sum_{r_2=1}^{R^{(k,l)}_2} \sum_{r_1=1}^{R^{(k,l)}_1}[\core{Q}^{(k,l)}]_{r_1,r_2,r_3} \matb{C}_{(3,r_3)}^{(k,l)}\otimes \matb{C}_{(2,r_2)}^{(k,l)}\otimes \matb{C}_{(1,r_1)}^{(k,l)},
\end{equation}
where  $\matb{C}_{(t,r_t)}^{(k,l)}\in\mathbb{R}^{n_t\times n_t}$ are matrices defined as
$$
\left[\matb{C}_{(t,r_t)}^{(k,l)}\right]_{i,j}:=\int_0^1   \delta^{(l,t)}(\widehat{b}_{i}(\eta_t) )\delta^{(k,t)}(\widehat{b}_{j}(\eta_t))\matb{Q}^{(k,l)}_{t}(:,r_t) \cdot {\vect{T}}_{p_t^{(k,l)}}(\eta_t) \d\eta_t.
$$
Note that $\widetilde{\matb{A}}$ is a Tucker matrix itself: by defining a single large block-diagonal core tensor  $\widetilde{\core{A}}\in\mathbb{R}^{R^A_{1}\times R^A_{2}\times R^A_{3}}$ with $R^A_{s}:=\sum_{k,l=1}^3 R_s^{(k,l)}$, having the tensors $\core{Q}^{(k,l)}$ for $k,l=1,2,3,$ as diagonal blocks, and factor matrices $\matb{A}_{(t,i)}\in\mathbb{R}^{n_t\times n_t}$ defined consequently from  $\matb{C}_{(t,r_t)}^{(k,l)}$, we get
\begin{equation*} 
\widetilde{\matb{A}} =  \sum_{r_3=1}^{R^A_{3}}\sum_{r_2=1}^{R^A_{2}} \sum_{r_1=1}^{R^A_{1}}[\widetilde{\core{A}}]_{r_1,r_2,r_3} \widetilde{\matb{A}}_{(3,r_3)} \otimes\widetilde{\matb{A}}_{(2,r_2)}\otimes \widetilde{\matb{A}}_{(1,r_1)}. 
\end{equation*}

Similarly, replacing $\omega(\vect{\eta})$ by its approximation \eqref{eq:approxW} in \eqref{eq:rhs_param} and associating  the index $i$   with the indexes of the univariate functions, i.e. $i\leftrightarrow (i_1,i_2,i_3)$,  we get
 \begin{align*} 
[\vect{f}]_{i} \approx[\widetilde{\vect{f}}]_{i}&:=  \int_{\widehat{\Omega}}   \left[\left({\vect{T}}_{q_3}(\eta_3)\matb{W}_3\otimes {\vect{T}}_{q_2}(\eta_2)\matb{W}_2\otimes {\vect{T}}_{q_1}(\eta_1)\matb{W}_1 \right) \vec(\core{W} )\right]
 \widehat{B}_i(\vect{\eta})  \d\widehat{\Omega}\\ 
&=  \sum_{r_3=1}^{R_3} \sum_{r_2=1}^{R_2} \sum_{r_1=1}^{R_1}[\core{F}]_{r_1,r_2,r_3}   \prod_{t=1}^3 \int_0^1 \widehat{b}_{t,i_t}(\eta_t) \matb{W}_t(:,r_t) \cdot{\vect{T}}_{q_t}(\eta_t)\d\eta_t,
 \end{align*}
 where $\matb{W}_t(:,r_t)\in\mathbb{R}^{m_t}$ denotes the $r_t$-th column of $\matb{S}_t$.
 Finally, we conclude that $\widetilde{\vect{f}}$ is a Tucker vector associated to the Tucker tensor
 \begin{equation*} 
\widetilde{\mathcal{F}}:= \widetilde{\core{F}}\times_3\widetilde{\factor{F}}_3\times_2 \widetilde{\factor{F}}_2\times_1\widetilde{\factor{F}}_1,
 \end{equation*}
 where 
 $$ \widetilde{\core{F}}:=\core{W}\quad \text{and}\quad
[ \widetilde{\factor{F}}_t]_{i,j}:=\int_0^1 \widehat{b}_{i}(\eta_t) \matb{W}_t(:,j) \cdot{\vect{T}}_{q_t}(\eta_t)\d\eta_t.
 $$ 

\paragraph{Computational cost}
  Suppose that a function $g$ can be approximated as
 $$g(\boldsymbol{\eta})\approx \widetilde{g}(\boldsymbol{\eta}) :=  \sum_{i_1=1}^{r}\sum_{i_2=1}^{r}\sum_{i_3=1}^{r}[\mathfrak{G}]_{i_1,i_2,i_3}u_{1,i_1}(\eta_1)u_{2,i_2}(\eta_2) u_{3,i_3}(\eta_3) ,$$
where $(r,r,r)$ is the multilinear rank,   $u_{m,j}(\eta_m)=\sum_{k=1}^{t} c_{k,j}T_{k}(\eta_m)$ and where $t$ denotes the polynomial degree. Then, the number of function evaluations performed by the \textsf{chebfun3f} function is $O(t r + r^3)$ \cite{Dolgov2021}. 
This number depends on the degree $t$, but in all our numerical tests we observed $t \ll n_m$. As a consequence, the computational cost to approximate the matrix and the right-hand side is negligible with respect to the overall cost of our approach.

\section{Approximation of eigenvalues and eigenvectors}
\label{sec:appendix}

In this appendix, we describe a way to approximate the eigenvector and eigenvalue matrices appearing in the Fast Diagonalization method (see Section \ref{sec:preconditioner} and Algorithm \ref{al:FD}), so that the corresponding matrix-vector products can be computed with almost-optimal complexity, while preserving the robustness with respect to $h$ and $p$ of the overall FD preconditioner. 
Here we limit the discussion to a description of the method, and postpone a throughout analysis to a forthcoming publication.


Let $\widehat{\S}_{h}^p$ be the space of splines on the interval $[0,1]$ of degree $p$ and regularity $p-1$ built from a  uniform knot vector. We denote with $\widehat{\mathcal{V}}_h$ the subspace of $\widehat{\S}_{h}^p$ whose functions vanish at the extremes where Dirichlet boundary condition is prescribed.

We start by discussing the case of Dirichlet boundary condition on both extremes 0 and 1. Let $n := \text{dim}(\widehat{\mathcal{V}}_{h}) = n_{el} + p - 2$, where $n_{el}$ denotes the number of knot spans. We assume $n_{el} > p $.

Let $\matb{M}$ and $\matb{K} \in \mathbb{R}^{n \times n}$ denote the mass and stiffness matrices for the univariate Poisson problem associated to $\widehat{\mathcal{V}}_{h}$. We consider the generalized eiegendecomposition
$$ \matb{K} \matb{U} = \matb{M} \matb{U} \boldsymbol{\Lambda}, \qquad \mbox{with} \qquad \matb{U}^T \matb{M} \matb{U} = \mathbf{I}_n,$$
where $\matb{U}$ is the ($\matb{M}-$orthogonal) eigenvector matrix and $\boldsymbol{\Lambda}$ is the diagonal matrix of eigenvalues.

When $p=1,2$, the matrices $\matb{U}$ and $\boldsymbol{\Lambda}$ are known explicitly \cite{Ekstrom2018}. Due to the peculiar structure of $\matb{U}$, computing a matrix-vector product with this matrix is equivalent to an application of a variant of Discrete Fourier Transform. This can be done efficiently thanks to the Fast Fourier Transform algorithm whose computational cost is $O(n \log n)$ \cite{VanLoan1992}. Note that no approximation of eigenvalues and eigenvectors is needed in this case.

For $p \geq 3$, we consider the following splitting of the spline space (note that the same splitting considered in \cite{Hofreither2017}) 
$$\widehat{\mathcal{V}}_{h}= \widehat{\mathcal{V}}_{h,1}\oplus \widehat{\mathcal{V}}_{h,2}, $$
where 
$$ \widehat{\mathcal{V}}_{h,1}:= \lbrace \phi \in \widehat{\mathcal{V}}_{h} \; \vert \; \phi^{(2k)}(0) = \phi^{(2k)}(1) = 0, \mbox{ for } k \in \mathbb{N} \mbox{ with } 2k \leq p-1  \rbrace $$
and $\widehat{\mathcal{V}}_{h,2}$ is the $\matb{M}-$orthogonal complement of $\widehat{\mathcal{V}}_{h,1}$ in $\widehat{\mathcal{V}}_{h}$.  We remark that $\widehat{\mathcal{V}}_{h,1}$  has been recently analyzed e.g. in  \cite{Deng2021,Hiemstra2021,Manni2022}. 
Note that
$$ n_{1} := \text{dim}(\widehat{\mathcal{V}}_{h,1}) = \left\{ \begin{array}{ll} n_{el} - 1 & \mbox{if } p \mbox{ is odd} \\ n_{el}  & \mbox{if } p \mbox{ is even} \end{array}\right. \quad \text{ and }\quad n_{2} := \text{dim}(\widehat{\mathcal{V}}_{h,2}) = \left\{ \begin{array}{ll} p - 1 & \mbox{if } p \mbox{ is odd} \\ p-2 & \mbox{if } p \mbox{ is even} \end{array} \right. .$$

We fix a basis for $\widehat{\mathcal{V}}_{h,1}$ that promotes sparsity. Specifically, if $\widehat{b}_1,\ldots, \widehat{b}_n$ denote the standard basis functions of $\widehat{\mathcal{V}}_{h}$, let $ \Phi = \left\{\widehat{b}_{n_{2} + 1}, \ldots, \widehat{b}_{n_{1}} \right\}$, and let $\Phi_{0}$ be a basis for the subspace
\begin{equation} \label{eq:Phi0} \mbox{span} \left\{\widehat{b}_1, \ldots, \widehat{b}_{n_{2}} \right\} \cap \left\{ \phi \in \widehat{\mathcal{V}}_{h}  \; \vert \; \phi^{(2k)}(0) = 0,  \mbox{ for } k \in \mathbb{N} \mbox{ with } 2k \leq p-1 \right\}. \end{equation}
Similarly, let $\Phi_{1}$ be a basis for the subspace
\begin{equation} \label{eq:Phi1} \mbox{span} \left\{ \widehat{b}_{n_{1} + 1}, \ldots, \widehat{b}_n \right\} \cap \left\{ \phi \in \widehat{\mathcal{V}}_{h} \; \vert \; \phi^{(2k)}(1) = 0, \mbox{ for } k \in \mathbb{N} \mbox{ with } 2k \leq p-1  \right\}. \end{equation}
Then we consider $ \Phi_0 \cup \Phi \cup \Phi_1$ as a basis for $\widehat{\mathcal{V}}_{h,1} $, and let $\matb{V}_{1} \in  \mathbb{R}^{n \times n_{2}}$ be the matrix whose columns represent these basis functions.
Note that the columns of $\matb{V}_{1}$ corresponding to functions in $\Phi_0$ and $\Phi_1$ are dense, while each column corresponding to a function in $\Phi$ has only one nonzero element.
Similarly, let $\matb{V}_{2} \in \mathbb{R}^{n \times n_{2}}$ be a matrix whose columns represent a basis for $\widehat{\mathcal{V}}_{h,2} $.

Since $n_{2}$ is small, we can afford to explicitly compute the eigendecomposition of the problem projected on $\widehat{\mathcal{V}}_{h,2} $, i.e. we compute  $\matb{U}_{2},\boldsymbol{\Lambda}_{2} \in \mathbb{R}^{n_{2} \times n_{2}}$ such that
$$ \left( \matb{V}_{2}^T \matb{K} \matb{V}_{2} \right) \matb{U}_{2} = \left( \matb{V}_{2}^T \matb{M} \matb{V}_{2} \right) \matb{U}_{2} \boldsymbol{\Lambda}_{2}. $$

On $\widehat{\mathcal{V}}_{h,1}$, on the other hand, we approximate the eigenvectors of the discrete operator using the interpolates of the eigenvectors of the analytic operator, which are sinusoidal functions.

More precisely, let $v_1,\ldots,v_{n_{1}}$ denote the basis functions of $\widehat{\mathcal{V}}_{h,1} $ defined above, and let $x_1, \ldots, x_{n_{1}}$ denote the interpolation points, defined as the breakpoints of the knot vector (except 0 and 1) for odd $p$, and as the midpoints of the knot spans for even $p$. We consider the collocation matrix $\matb{C} \in \mathbb{R}^{n_{1} \times n_{1}}$ with entries 
$$ [\matb{C}]_{i,j} = v_i(x_j) \qquad \text{for }i,j = 1,\ldots,n_{1}$$
and the matrix $\matb{U}_{dft} \in \mathbb{R}^{n_{1} \times n_{1}} $ whose entries are sinusoidal functions (normalized in the $L^2$ norm) evaluated at the interpolation points, i.e.
$$ \left[ \matb{U}_{dft} \right]_{i,j} = \sqrt{2} \sin(j \pi x_i ) \qquad \text{for } i,j = 1,\ldots,n_{1}.$$

Then on $\widehat{\mathcal{V}}_{h,1}$ the matrix of eigenvector is approximated with 
$ \matb{U}_{1} := \sqrt{2} \matb{C}^{-1} \matb{U}_{dft}$. 

The approximated eigenvalue matrix $\boldsymbol{\Lambda}_{1}$, on the other hand, is constructed using the eigenvalues of the analytic operator, i.e.

$$ \left[\boldsymbol{\Lambda}_{1} \right]_{j,j} = (j \pi)^2 \qquad \text{for } j = 1,\ldots,n_{1}. $$

In conclusion, the approximation to the full eigenvector matrix $\matb{U}$ is taken as:
$$ \widetilde{\matb{U}} := \begin{bmatrix} \matb{V}_{1} \matb{U}_{1} \; \vline \; \matb{V}_{2} \matb{U}_{2} \end{bmatrix}$$

Similarly, the approximation of the eigenvalue matrix $\Lambda$ is taken as:
$$ \widetilde{\matb{\Lambda}} := \begin{bmatrix} \boldsymbol{\Lambda}_{1} & 0 \\ 0 & \boldsymbol{\Lambda}_{2} \end{bmatrix}$$

We emphasize that matrix-vector products with $\widetilde{\matb{U}}$ can be computed with $O( n (\log (n) + p))$ complexity.
The key point is that the matrix $\matb{U}_{dft}$ does not need to be computed and stored, and its action on a vector can be computed exploiting the FFT, which costs $n \log (n)$ operations. The action of $\matb{C}^{-1}$ requires the solution of a banded linear system with $p+1$ nonzero entries per row, which costs $O(np)$ operations. This also the cost of matrix-vector products with $\matb{V}_{1}$ and $ \matb{V}_{2} \matb{U}_{2} \in \mathbb{R}^{n \times n_{2}}$.  

If Neumann boundary condition is present, and assuming $p \geq 2$ (for $p=1$ eigenvalues and eigenvectors are explicitly known, see \cite[Section 4.5]{VanLoan1992}), the definition of $\widehat{\mathcal{V}}_{h,1}$ has to be modified. Here the functions of $\widehat{\mathcal{V}}_{h,1}$ have vanishing \textit{odd} derivatives at the extremes with Neumann boundary condition (and vanishing even derivatives at the extremes with Dirichlet boundary condition).
$\widehat{\mathcal{V}}_{h,2}$ is still defined as the $\matb{M}-$orthogonal complement of $\widehat{\mathcal{V}}_{h,1}$, and the construction of $\matb{V}_2$, $\matb{U}_2$ and $\matb{V}_1$ proceeds analogously as before (note that the definition of the spaces \eqref{eq:Phi0} and \eqref{eq:Phi1} has to be modified according to $\widehat{\mathcal{V}}_{h,1}$).

The approximated eigenvectors on $\widehat{\mathcal{V}}_{h,1}$ are still sinusoidal functions satisfying the derivative restrictions embedded into $\widehat{\mathcal{V}}_{h,1}$ at the boundary. Precisely we have
$$ \left[ \matb{U}_{dft} \right]_{i,j} = \sin \left( \left( j - \frac{k_0}{2} - \frac{k_1}{2} \right) \pi x_i + \frac{k_0 \pi}{2} \right) \qquad \text{for } i,j = 1,\ldots,n_{1}, $$
where
$$ k_{0} = \left\{ \begin{array}{ll} 1 & \mbox{if a Neumann b.c. is prescribed at } 0 \\ 0  & \mbox{otherwise} \end{array}\right., \qquad k_{1} = \left\{ \begin{array}{ll} 1 & \mbox{if a Neumann b.c. is prescribed at } 1 \\ 0 & \mbox{otherwise} \end{array} \right.,$$
and where the interpolation points $x_i$, $i=1,\ldots,n_1$, are still the midpoints of the knot spans for even $p$ or the breakpoints of the knot vector (excluding any extreme associated with Dirichlet boundary condition) for odd $p$.
The approximated eigenvalues are then 
$$ \left[\boldsymbol{\Lambda}_{1} \right]_{j,j} = \left( j - \frac{k_0}{2} - \frac{k_1}{2} \right)^2 \pi^2 \qquad \text{for } j = 1,\ldots,n_{1}. $$
The interpolation matrix $\matb{C}$ and the final matrices $\widetilde{\matb{U}}$ and $\widetilde{\matb{\Lambda}}$ are defined as before.
The analysis of the computational cost is also analogous as the previous case. Note in particular that the FFT can still be exploited to compute matrix-vector products with $\matb{U}_{dft}$ (see \cite{VanLoan1992} for details).

\bibliographystyle{plain}
\bibliography{biblio_tensor_solver}

\begin{thebibliography}{10}

\bibitem{Beirao2014}
L.~Beir{\~a}o~da Veiga, A.~Buffa, G.~Sangalli, and R.~V\'{a}zquez.
\newblock Mathematical analysis of variational isogeometric methods.
\newblock {\em Acta Numerica}, 23:157--287, 2014.

\bibitem{Beirao2017}
L.~Beir{\~a}o~da Veiga, L.~F. Pavarino, S.~Scacchi, O.~B. Widlund, and
  S.~Zampini.
\newblock Adaptive selection of primal constraints for isogeometric bddc deluxe
  preconditioners.
\newblock {\em SIAM Journal on Scientific Computing}, 39(1):A281--A302, 2017.

\bibitem{Bosy2020}
M.~Bosy, M.~Montardini, G.~Sangalli, and M.~Tani.
\newblock A domain decomposition method for isogeometric multi-patch problems
  with inexact local solvers.
\newblock {\em Computers \& Mathematics with Applications}, 80(11):2604--2621,
  2020.

\bibitem{Braess2005}
D.~Braess and W.~Hackbusch.
\newblock Approximation of $1/x$ by exponential sums in $[1, \infty)$.
\newblock {\em IMA Journal of Numerical Analysis}, 25(4):685--697, 2005.

\bibitem{Bunger2020}
A.~B{\"u}nger, S.~Dolgov, and M.~Stoll.
\newblock A low-rank tensor method for {PDE}-constrained optimization with
  isogeometric analysis.
\newblock {\em SIAM Journal on Scientific Computing}, 42(1):A140--A161, 2020.

\bibitem{Cottrell2009}
J.~A. Cottrell, T.~J.~R. Hughes, and Y.~Bazilevs.
\newblock {\em Isogeometric analysis: toward integration of {CAD} and {FEA}}.
\newblock John Wiley \& Sons, 2009.

\bibitem{Cottrell2007}
J.~A. Cottrell, T.~J.~R. Hughes, and A.~Reali.
\newblock Studies of refinement and continuity in isogeometric structural
  analysis.
\newblock {\em Computer Methods in Applied Mechanics and Engineering},
  196(41):4160--4183, 2007.

\bibitem{DeBoor2001}
C.~De~Boor.
\newblock {\em {A practical guide to splines (revised edition)}}.
\newblock Applied {M}athematical {S}ciences. Springer, 2001.

\bibitem{De2000}
L.~De~Lathauwer, B.~De~Moor, and J.~Vandewalle.
\newblock A multilinear singular value decomposition.
\newblock {\em SIAM journal on Matrix Analysis and Applications},
  21(4):1253--1278, 2000.

\bibitem{Deng2021}
Q.~Deng and V.~M. Calo.
\newblock A boundary penalization technique to remove outliers from
  isogeometric analysis on tensor-product meshes.
\newblock {\em Computer Methods in Applied Mechanics and Engineering},
  383:113907, 2021.

\bibitem{Dolgov2021}
S.~Dolgov, D.~Kressner, and C.~Strossner.
\newblock Functional {T}ucker approximation using {C}hebyshev interpolation.
\newblock {\em SIAM Journal on Scientific Computing}, 43(3):A2190--A2210, 2021.

\bibitem{Donatelli2015}
M.~Donatelli, C.~Garoni, C.~Manni, S.~Serra-Capizzano, and H.~Speleers.
\newblock Robust and optimal multi-iterative techniques for {I}g{A} {G}alerkin
  linear systems.
\newblock {\em Computer Methods in Applied Mechanics and Engineering},
  284:230--264, 2015.

\bibitem{Driscoll2014}
T.~A Driscoll, N.~Hale, and L.~N. Trefethen.
\newblock {\em Chebfun Guide}.
\newblock Pafnuty Publications, 2014.

\bibitem{Ekstrom2018}
S.~Ekstr{\"o}m, I.~Furci, C.~Garoni, C.~Manni, S.~Serra-Capizzano, and
  H.~Speleers.
\newblock Are the eigenvalues of the {B}-spline isogeometric analysis
  approximation of- {$\Delta$} u= $\lambda$ u known in almost closed form?
\newblock {\em Numerical Linear Algebra with Applications}, 25(5):e2198, 2018.

\bibitem{Garcia2018}
D.~Garcia, D.~Pardo, L.~Dalcin, and V.~M. Calo.
\newblock Refined {I}sogeometric {A}nalysis for a preconditioned conjugate
  gradient solver.
\newblock {\em Computer Methods in Applied Mechanics and Engineering},
  335:490--509, 2018.

\bibitem{Georgieva2019}
I.~Georgieva and C.~Hofreither.
\newblock Greedy low-rank approximation in {T}ucker format of solutions of
  tensor linear systems.
\newblock {\em Journal of Computational and Applied Mathematics}, 358:206--220,
  2019.

\bibitem{Grasedyck2013}
L.~Grasedyck, D.~Kressner, and C.~Tobler.
\newblock A literature survey of low-rank tensor approximation techniques.
\newblock {\em GAMM-Mitteilungen}, 36(1):53--78, 2013.

\bibitem{Hackbusch2019}
W.~Hackbusch.
\newblock Computation of best {${L}^{\infty}$} exponential sums for $1/x$ by
  {R}emez’algorithm.
\newblock {\em Computing and Visualization in Science}, 20(1):1--11, 2019.

\bibitem{Hiemstra2021}
R.~R Hiemstra, T.~J.~R. Hughes, A.~Reali, and D.~Schillinger.
\newblock Removal of spurious outlier frequencies and modes from isogeometric
  discretizations of second-and fourth-order problems in one, two, and three
  dimensions.
\newblock {\em Computer Methods in Applied Mechanics and Engineering},
  387:114115, 2021.

\bibitem{Hofreither2018}
C.~Hofreither.
\newblock A black-box low-rank approximation algorithm for fast matrix assembly
  in isogeometric analysis.
\newblock {\em Computer Methods in Applied Mechanics and Engineering},
  333:311--330, 2018.

\bibitem{Hofreither2017}
C.~Hofreither and S.~Takacs.
\newblock Robust multigrid for isogeometric analysis based on stable splittings
  of spline spaces.
\newblock {\em SIAM Journal on Numerical Analysis}, 55(4):2004--2024, 2017.

\bibitem{Hughes2005}
T.~J.~R. Hughes, J.~A. Cottrell, and Y.~Bazilevs.
\newblock Isogeometric analysis: {CAD}, finite elements, {NURBS}, exact
  geometry and mesh refinement.
\newblock {\em Computer methods in applied mechanics and engineering},
  194(39-41):4135--4195, 2005.

\bibitem{HST_at_kunoth2018isogeometric}
T.~J.~R. Hughes, G.~Sangalli, and M.~Tani.
\newblock Isogeometric analysis: Mathematical and implementational aspects,
  with applications.
\newblock {\em Splines and PDEs: From Approximation Theory to Numerical Linear
  Algebra: Cetraro, Italy 2017}, pages 237--315, 2018.

\bibitem{Juttler2017}
B.~J{\"u}ttler and D.~Mokri{\v{s}}.
\newblock Low rank interpolation of boundary spline curves.
\newblock {\em Computer Aided Geometric Design}, 55:48--68, 2017.

\bibitem{Kolda2009}
T.~G. Kolda and B.~W. Bader.
\newblock Tensor decompositions and applications.
\newblock {\em SIAM review}, 51(3):455--500, 2009.

\bibitem{Kressner2017}
D.~Kressner and L.~Perisa.
\newblock Recompression of hadamard products of tensors in tucker format.
\newblock {\em SIAM Journal on Scientific Computing}, 39(5):A1879--A1902, 2017.

\bibitem{Kressner2010}
D.~Kressner and C.~Tobler.
\newblock Krylov subspace methods for linear systems with tensor product
  structure.
\newblock {\em SIAM Journal on Matrix Analysis and Applications},
  31(4):1688--1714, 2010.

\bibitem{Kressner2011}
D.~Kressner and C.~Tobler.
\newblock Low-rank tensor {K}rylov subspace methods for parametrized linear
  systems.
\newblock {\em SIAM Journal on Matrix Analysis and Applications},
  32(4):1288--1316, 2011.

\bibitem{Loli2022}
G.~Loli, G.~Sangalli, and M.~Tani.
\newblock Easy and efficient preconditioning of the isogeometric mass matrix.
\newblock {\em Computers \& Mathematics with Applications}, 116:245--264, 2022.

\bibitem{Lynch1964}
R.~E. Lynch, J.~R. Rice, and D.~H. Thomas.
\newblock Direct solution of partial difference equations by tensor product
  methods.
\newblock {\em Numerische Mathematik}, 6(1):185--199, 1964.

\bibitem{Manni2022}
C.~Manni, E.~Sande, and H.~Speleers.
\newblock Application of optimal spline subspaces for the removal of spurious
  outliers in isogeometric discretizations.
\newblock {\em Computer Methods in Applied Mechanics and Engineering},
  389:114260, 2022.

\bibitem{Mantzaflaris2017}
A.~Mantzaflaris, B.~J{\"u}ttler, B.~N. Khoromskij, and U.~Langer.
\newblock Low rank tensor methods in {G}alerkin-based isogeometric analysis.
\newblock {\em Comput. Methods Appl. Mech. Engrg.}, 316:1062--1085, 2017.

\bibitem{Matthies2012}
H.~G. Matthies and E.~Zander.
\newblock Solving stochastic systems with low-rank tensor compression.
\newblock {\em Linear Algebra and its Applications}, 436(10):3819--3838, 2012.

\bibitem{Niederreiter1992}
H.~Niederreiter.
\newblock {\em Random number generation and quasi-Monte Carlo methods}.
\newblock SIAM, 1992.

\bibitem{Oseledets2009}
I.~V. Oseledets, D.~V. Savostyanov, and E.~E. Tyrtyshnikov.
\newblock Linear algebra for tensor problems.
\newblock {\em Computing}, 85(3):169--188, 2009.

\bibitem{Palitta2021}
D.~Palitta and P.~K{\"u}rschner.
\newblock On the convergence of {K}rylov methods with low-rank truncations.
\newblock {\em Numerical Algorithms}, 88(3):1383--1417, 2021.

\bibitem{Pan2019}
M.~Pan and F.~Chen.
\newblock Low-rank parameterization of volumetric domains for isogeometric
  analysis.
\newblock {\em Computer-Aided Design}, 114:82--90, 2019.

\bibitem{Sande2020}
E.~Sande, C.~Manni, and H.~Speleers.
\newblock Explicit error estimates for spline approximation of arbitrary
  smoothness in isogeometric analysis.
\newblock {\em Numerische Mathematik}, 144(4):889--929, 2020.

\bibitem{Sangalli2016}
G.~Sangalli and M.~Tani.
\newblock Isogeometric preconditioners based on fast solvers for the
  {S}ylvester equation.
\newblock {\em SIAM Journal on Scientific Computing}, 38(6):A3644--A3671, 2016.

\bibitem{Simoncini2022}
V.~Simoncini and Y.~Hao.
\newblock Analysis of the truncated conjugate gradient method for linear matrix
  equations.
\newblock {\em hal-03579267}, 2022.

\bibitem{Simoncini2003}
V.~Simoncini and D.~B. Szyld.
\newblock Theory of inexact {K}rylov subspace methods and applications to
  scientific computing.
\newblock {\em SIAM Journal on Scientific Computing}, 25(2):454--477, 2003.

\bibitem{Tielen2020}
R.~Tielen, M.~M{\"o}ller, D.~G{\"o}ddeke, and C.~Vuik.
\newblock $p$-multigrid methods and their comparison to $h$-multigrid methods
  within {I}sogeometric {A}nalysis.
\newblock {\em Computer Methods in Applied Mechanics and Engineering},
  372:113347, 2020.

\bibitem{Tucker1966}
L.~R. Tucker.
\newblock Some mathematical notes on three-mode factor analysis.
\newblock {\em Psychometrika}, 31(3):279--311, 1966.

\bibitem{VanLoan1992}
C.~Van~Loan.
\newblock {\em Computational frameworks for the {F}ast {F}ourier {T}ransform}.
\newblock SIAM, 1992.

\bibitem{Vannieuwenhoven2012}
N.~Vannieuwenhoven, R.~Vandebril, and K.~Meerbergen.
\newblock A {N}ew {T}runcation {S}trategy for the {H}igher-{O}rder {S}ingular
  {V}alue {D}ecomposition.
\newblock {\em SIAM Journal on Scientific Computing}, 34(2):A1027--A1052, 2012.

\bibitem{Vazquez2016}
R.~V{\'a}zquez.
\newblock A new design for the implementation of isogeometric analysis in
  {O}ctave and {M}atlab: {G}eo{PDE}s 3.0.
\newblock {\em Computers \& Mathematics with Applications}, 72(3):523--554,
  2016.

\bibitem{Sorber2014}
N.~Vervliet, O.~Debals, L.~Sorber, M.~Van~Barel, and L.~De~Lathauwer.
\newblock Tensorlab v3.0.
\newblock {\em URL: www.tensorlab.net}, 2016.

\bibitem{Zander2013}
E.~K. Zander.
\newblock {\em Tensor approximation methods for stochastic problems}.
\newblock PhD thesis, Braunschweig, Institut f{\"u}r Wissenschaftliches
  Rechnen, 2013.

\end{thebibliography}

\end{document}